\documentclass[11pt,a4paper]{amsart}

\usepackage[utf8x]{inputenc}
\usepackage{ucs}
\usepackage{amsmath}
\usepackage{amsfonts}
\usepackage{amssymb}
\usepackage{graphicx}
\usepackage{geometry}
\usepackage{color}
\usepackage[active]{srcltx}
\usepackage[all]{xy}

\allowdisplaybreaks

\newtheorem{lema}{Lemma}[section]
\newtheorem{defi}[lema]{Definition}

\newtheorem{teor}[lema]{Theorem}
\newtheorem{cor}[lema]{Corollary}
\newtheorem{prop}[lema]{Proposition}

\newtheorem{obs}[lema]{Remark}
\newtheorem{rmks}[lema]{Remarks}

\newcounter{maint}

\newtheorem{mteor}[maint]{Theorem}

\newcommand{\Ss}{{\mathcal{S}}}
\newcommand{\Kk}{{\mathcal{K}}}
\newcommand{\Aa}{{\mathcal{A}}_4^{\prime\prime}}
\newcommand{\A}{\mathfrak{A}}
\newcommand{\com}{\Delta}
\newcommand{\eps}{\varepsilon}
\newcommand{\gr}{\operatorname{gr}}
\newcommand\co{\operatorname{co}}
\newcommand\op{\operatorname{op}}
\newcommand\cop{\operatorname{cop}}
\newcommand\id{\operatorname{id}}
\newcommand\ad{\operatorname{ad}}
\newcommand\End{\operatorname{End}}

\newcommand\Soc{\operatorname{Soc}}
\newcommand{\Z}{{\mathbb Z}}
\newcommand{\N}{{\mathbb N}}

\newcommand{\ydho}{{}^{H_0}_{H_0}\mathcal{YD}}
\newcommand{\ydk}{{}^{\Kk}_{\Kk}\mathcal{YD}}
\newcommand{\yda}{{}^{\mathcal{A}}_{\mathcal{A}}\mathcal{YD}}
\newcommand{\ydra}{{}^{A}_{A}\mathcal{YD}}
\newcommand{\ydh}{{}^H_H\mathcal{YD}}
\newcommand{\ydrho}{{}^{H_{[0]}}_{H_{[0]}}\mathcal{YD}}
\newcommand\toba{{\mathfrak B }}
\newcommand\cal{\mathcal }
\newcommand\Hom{\operatorname{Hom}}
\newcommand\Alg{\operatorname{Alg}}

\newcommand\ExtQ{\operatorname{ExtQ}}
\newcommand\Ext{\operatorname{Ext}}
\newcommand\Top{\operatorname{Top}}

\usepackage{soul}
\def\pf{\begin{proof}}
\def\epf{\end{proof}}
\def\ot{\otimes}

\begin{document}

\title[On Hopf algebras over quantum subgroups]
{On Hopf algebras over quantum subgroups}

\author[G. A. Garc\'ia]
{Gast\'on Andr\'es Garc\'ia}

\author[J. M. J. Giraldi]{Jo\~ao Matheus Jury Giraldi}

\thanks{
This work was partially supported by
ANPCyT-Foncyt, CONICET, SeCyT (UNLP), CNPq (Brazil)}


\address{\noindent G. A. G. : Departamento de Matem\'atica, Facultad de Ciencias Exactas,
Universidad Nacional de La Plata. CONICET. Casilla de Correo 172, (1900)
La Plata, Argentina.}

\address{\noindent J. M. J. G. : Instituto de Matem\'atica, 
Universidade Federal do Rio Grande do Sul,
Rio Grande do Sul, Brazil}

\email{ggarcia@mate.unlp.edu.ar,
joaomjg@gmail.com}

\begin{abstract}
Using the standard filtration associated with a generalized lifting method,
we determine all finite-dimensional Hopf algebras over an algebraically closed field of 
characteristic zero whose coradical generates a Hopf subalgebra isomorphic to 
the smallest non-pointed non-cosemisimple Hopf algebra $\Kk$ and the corresponding 
infinitesimal module is an indecomposable object in $\ydk$ (we assume that the diagrams are Nichols algebras). 
As a byproduct, we 
obtain new Nichols algebras of dimension 8 and new Hopf algebras of dimension 64.
\end{abstract}

\maketitle

\section*{Introduction}
Let $\Bbbk$ be an algebraically closed field of characteristic zero.
The problem of classifying all Hopf algebras over $\Bbbk$ of a given dimension was posed by Kaplansky in 1975 
\cite{K}. 
Some progress has been made but, in general, it is
a difficult question. One of the few general techniques is the so-called \textit{Lifting Method} \cite{AS},
under the assumption that the coradical is a subalgebra, \textit{i.e.}, the Hopf algebra
has the Chevalley Property. More recently, 
Andruskiewitsch and Cuadra \cite{AC}
proposed 
to extend this technique by considering the subalgebra generated by
the coradical and the related wedge filtration.
It turns out that 
this
filtration is a Hopf algebra filtration,
provided that the antipode is injective, what is true in the finite-dimensional context.

We describe the lifting method briefly.
Let $H$ be a Hopf algebra over $\Bbbk$. Recall that the coradical filtration $\{H_{ n } \}_{n\geq 0} $ of $H$ 
is defined recursively by 
 \begin{itemize}
 \item the coradical $H_{0}$, which is the sum of all simple subcoalgebras, and
 \item $H_{n} = \bigwedge^{n+1}H_{0} = \{h\in H:\ \com(h) \in H\ot H_{0} + H_{n-1}\ot H\}$.
\end{itemize}
This filtration corresponds
to the filtration of $H^{*}$ given by the powers of the Jacobson radical.  It is 
always a coalgebra filtration and if $H_{ 0}$ is a
Hopf subalgebra, then it is indeed a Hopf algebra filtration; in particular, 
its associated graded object $\gr H = \bigoplus_{ n\geq 0} H_{ n} /H_{ n-1}$ is a graded Hopf
algebra, where  $H_{-1} = 0$. Let $\pi : \gr H \to H_{ 0}$ be the homogeneous projection.
It turns out that $\gr H \simeq R\# H_{0}$ as Hopf algebras, where
$R = (\gr H)^{\co \pi} = \{h\in H:\ (\id\ot \pi)\com(h)= h\ot 1\}$ is the algebra of coinvariants and 
$\#$ stands for the Radford-Majid biproduct or \textit{bosonization} of $R$ with $H_{0}$. 
The algebra $R$ is not a usual Hopf algebra, but a graded connected Hopf algebra in 
the category $\ydho$ of left Yetter-Drinfeld modules over $H_{0}$. The subalgebra generated by the elements of degree one is the \textit{Nichols algebra} $\toba(V)$  of $V = R(1)$; here 
$V$ is a braided vector space called the \textit{infinitesimal braiding}.

Let us fix a finite-dimensional cosemisimple Hopf algebra $A$.
The lifting method then consists of the description of all finite-dimensional Nichols algebras $ \toba(V) \in \ydra $, the determination of all possible deformations of the bosonization $\toba(V)\#A$, and the proof that all Hopf algebras $H$ with $H_{0}=A$ satisfy that $\gr H \simeq \toba(V)\#A$. 

The main idea in \cite{AC} is to 
replace the coradical filtration by a more general but adequate filtration: the 
\textbf{standard filtration} $\{H_{ [n]} \}_{ n\geq 0} $, which is defined recursively by
\begin{itemize}
 \item  the subalgebra $H_{[0]}$ of $H$ generated by $H_{ 0}$, called the \textit{Hopf coradical}, and  
 \item $H_{[n]} = \bigwedge^{n+1}H_{[0]}$.
\end{itemize}
If the coradical $H_{0}$ is a Hopf subalgebra, then  
$H_{ [0]} = H_{ 0}$ and the coradical filtration coincides with the standard one. 

Let $A$ be a Hopf algebra generated by its coradical. We will say that $H$ is a \textit{Hopf algebra over} $A$ if
$H_{[0]} \simeq A$ as Hopf algebras.

Assume that the antipode $\Ss$ of $H$ is injective. Then by \cite[Lemma 1.1]{AC}, it holds
that $H_{[0]}$ is a Hopf subalgebra of $H$, $H_{n}\subseteq H_{[n]}$ and 
$\{H_{ [n]} \}_{ n\geq 0} $ is a Hopf algebra filtration of $H$. In particular, 
the graded algebra
$\gr H = \bigoplus_{ n\geq 0} H_{ [n]} /H_{ [n-1]}$ with $H_{[-1]}=0$ is a Hopf algebra 
associated with the standard filtration. 
Write $\pi : \gr H \to H_{ [0]}$ for the homogeneous projection. Then, as before, it splits the inclusion
of $H_{ [0]}$ in $\gr H$, the \textit{diagram} $R = (\gr H)^{\co\pi}$ is a Hopf algebra in the
category $\ydrho$ of Yetter-Drinfeld modules over $H_{[0]}$ and 
$\gr H \simeq R\#H_{[0]}$ as Hopf algebras. It turns out that 
$R = \bigoplus_{n\geq 0}R(n)$ is also graded and connected. 
We call again the linear space $R(1)$ consisting of elements of degree one,
the \textit{infinitesimal braiding}. 

The procedure to describe explicitly any Hopf algebra as above 
defines a proposal for the classification
of general finite-dimensional Hopf algebras over a fixed Hopf subalgebra $A$
which is generated by a cosemisimple coalgebra. The 
main steps are the following:
\begin{enumerate}
\item[$(a)$] determine all Yetter-Drinfeld modules $V$ in $\ydra$ such that the
Nichols algebra $\toba(V)$ is finite-dimensional,

\item[$(b)$] for such $V$, compute all Hopf algebras $L$ such
that $\gr L\simeq \toba(V)\# A$. 
We call $L$ a \emph{lifting}
of $\toba(V)$ over $A$.

\item[$(c)$] Prove that any finite-dimensional 
Hopf algebra over $A$
is generated by the first term of the standard filtration. 
\end{enumerate}

In this paper, we study these questions ($a$) and ($b$) in the case that $A =\Kk$ is the smallest 
Hopf algebra whose coradical is not a subalgebra. It is 
an 8-dimensional Hopf algebra whose dual is a pointed Hopf algebra.
The dual Hopf algebra $A^{*}$ was first introduced by Radford \cite{R1}, who addressed
the problem of finding a Hopf algebra whose
Jacobson radical is not a Hopf ideal. 

Let $\xi$ be a primitive $4$-th root of 1.
As an algebra, $\Kk$ is generated by the elements $a,b,c,d$
satisfying the following relations:
\begin{align}\label{eq:defining_relations_K}
ab & = \xi ba, & ac & = \xi ca, & 0& =cb=bc, &  cd & = \xi dc,&
bd & = \xi db,\\
\notag ad & = da, &  ad &=1,
& 0& =b^{2}=c^{2}, & a^{2}c & = b, &  a^{4} & =1.
\end{align}
The coalgebra structure and its antipode are determined by 
\begin{align}
\notag \Delta(a) &= a\otimes a + b\otimes c, &\Delta(b) &= a\otimes b + b\otimes d, & \eps(a)& =1,& \eps(b) & =0,\\\label{eq:defining_coproduct_antipode_K}
 \Delta(c) & = c\otimes a + d\otimes c, & \Delta(d) & = c\otimes b + d\otimes d, & \eps(c)& =0,& \eps(d) & =1,\\
\notag  \Ss(a) & = d, & \Ss(b) &= \xi b, & \Ss(c) &= -\xi c, & \Ss(d) &= a.
\end{align}
See Section \ref{sec:qsubgroupK} for more details.

In order to determine finite-dimensional Hopf algebras over $\Kk$, we first compute 
the Drinfeld double $D:=D(\Kk^{\cop}) $ of $\Kk^{\cop}$ and describe  
 the simple left $D$-modules,
their projective covers and some indecomposable  left $D$-modules. In fact,
we prove in Theorem \ref{thm:simpleD-mod} that 
there are sixteen simple left  $ D $-modules pairwise non-isomorphic: four
$1$-dimensional ones and twelve $2$-dimensional ones. The former correspond to characters on $\Z_{4}$ and the latter are parametrized by the
set $\Lambda = \{(i,j) \in \Z_{4}\times \Z_{4}|\  
2i\neq j\}$.
We compute the separation diagram of $D$ and show that $D$ is of tame
representation type.

Using that the braided monoidal categories $_{D}\mathcal{M}$ and $\ydk$ are equivalent, we then translate
the description above to simple and indecomposable modules in $\ydk$.
Then, using the 
description of the braiding in $\ydk$, we obtain our first main result, see Section \ref{sec:categoryyd}
for definitions.

\begin{mteor}\label{thm:Nicholskdm}
Let $\toba(V)$ be a finite-dimensional Nichols algebra over an indecomposable object $V$ in $\ydk$. Then 
$ V$ is simple and isomorphic either to $\Bbbk_{\chi}$, $\Bbbk_{\chi^{3}}$, 
$V_{2,1}$, $V_{2,3}$, $V_{3,1}$ or $V_{3,3}$. 
\end{mteor}

It turns out that $\toba(\Bbbk_{\chi^{\ell}}) \simeq \bigwedge \Bbbk_{\chi^{\ell}}$ is an exterior 
algebra for $\ell = 1,3$ with $\dim \toba(\Bbbk_{\chi^{\ell}}) = 2$ and $\toba(V)$ is an $8$-dimensional
algebra for $V= V_{2,1}$, $V_{2,3}$, $V_{3,1}$ and $V_{3,3}$.
It is possible to check that these braidings are triangular \cite{ufer}. 
These 8-dimensional examples are new examples of finite-dimensional Nichols algebras. They 
are isomorphic to quantum linear spaces as algebras, but not as coalgebras since the
braiding differs; in our case, the braiding is not of diagonal type, see the Appendix of the first
\texttt{arXiv} version of this paper. 

As the study of Nichols algebras over semisimple modules is a hard problem that
demands different techniques to be applied, we focus on the description of 
Hopf algebras over $\Kk$ such that their infinitesimal braiding is indecomposable, 
\textit{i.e.}, the liftings of the Nichols algebras in Theorem \ref{thm:Nicholskdm}. 
Thus we define two Hopf algebras $\A_{3,1}(\mu)$ and $\A_{3,3}(\mu)$ depending
on a parameter $\mu \in \Bbbk$ and prove our second main result, see Section \ref{sec:HaK} for 
definitions.

\begin{mteor}\label{mthm:liftingsoverK}
 Let $H$ be a finite-dimensional Hopf algebra over $\Kk$ such that its 
 infinitesimal braiding is an indecomposable module $V$ in $\ydk$.
Assume that the diagram $ R $ is a Nichols algebra. 
Then $V$ is simple and $H$ is isomorphic either to
 \begin{enumerate}
  \item[$(i)$] $(\bigwedge \Bbbk_{\chi^{\ell}})\# \Kk$ with $\ell = 1,3$;
  \item[$(ii)$] $\toba(V_{2,1})\# \Kk$;
  \item[$(iii)$] $\toba(V_{2,3})\# \Kk$;
  \item[$(iv)$] $\A_{3,1}(\mu)$ for some $\mu \in \Bbbk$;
  \item[$(v)$] $\A_{3,3}(\mu)$ for some $\mu \in \Bbbk$.
 \end{enumerate}
\end{mteor}

The Hopf algebras $(\bigwedge \Bbbk_{\chi^{\ell}})\# \Kk$ with $\ell = 1,3$ have dimension 
 $16$ and are duals of pointed Hopf algebras. They have already appeared in \cite{B}.
 The Hopf algebras $\toba(V_{2,1})\# \Kk$ and $\toba(V_{2,3})\# \Kk$ are
 dual of pointed Hopf algebras of dimension $64$.
 The Hopf algebras $\A_{3,1}(\mu)$ and $\A_{3,3}(\mu)$ are non-pointed with non-pointed duals. 
 To the best of the authors knowledge, they constitute new examples of Hopf 
 algebras of dimension $64$. 

 The paper is organized as follows. In Section \ref{sec:prelim} we recall some  
 invariants associated with a Hopf algebra, define Yetter-Drinfeld modules, Nichols algebras and the 
 Drinfeld double of a finite-dimensional Hopf algebra. We also recall
 the relation between Hopf algebras with a projection and bosonizations. In Section \ref{sec:qsubgroupK}
we describe the structure of $\Kk$ and 
give the presentation of the double $D= D(\Kk^{\cop})$ by generators and relations. We also 
determine the simple left $D$-modules,
their projective covers and some indecomposable  left $D$-modules. 
We compute the Ext-Quiver of $D$ and show that $D$
is of tame representation type.

Then, using 
the equivalence $_{D}\mathcal{M} \simeq \ydk$, we determine in Section \ref{sec:categoryyd}
the corresponding objects of the latter and describe their braidings. 
In Section \ref{sec:NicholsalgK} we show that if $\toba(V)$ is a 
finite-dimensional Nichols algebra in $\ydk$, then $V$ is necessarily 
a semisimple object and prove Theorem \ref{thm:Nicholskdm} by describing
first the Nichols algebra of the simple modules. Finally, in Section \ref{sec:HaK} we prove Theorem \ref{mthm:liftingsoverK}.

\section{Preliminaries}\label{sec:prelim}

\subsection{Conventions}
We work over an algebraically closed field $\Bbbk$ of characteristic
zero and with Hopf algebras which have bijective antipode. Our references for Hopf algebra theory are \cite{Mo} and
\cite{R2}.

For a Hopf algebra $H$ over $\Bbbk$, the comultiplication,
counit and antipode are denoted by $\com$, $\eps$
and $\Ss$, respectively. 
Comultiplication and coactions are written using the 
Sweedler  notation with summation sign suppressed, \textit{e.g.},
$\com(h)=h_{(1)}\ot
h_{(2)}$ for $h\in H$. 
A Hopf algebra in a braided monoidal category is called a 
\textit{braided} Hopf algebra. 
We denote by $_{H}\mathcal{M}$ the category of finite-dimensional
left $H$-modules.

The set $G(H)=\{h\in H\setminus\{0\}:\ \com(h) = h\ot h\}$
denotes the group
of \emph{group-like elements}.
The \emph{coradical}
$H_0$ of $H$ is the sum of all simple
subcoalgebras of $H$; in particular, $\Bbbk G(H)\subseteq
H_0$. The subalgebra $H_{[0]}$ generated by $H_{0}$
is a Hopf subalgebra which is called the \emph{Hopf coradical}.
For $h,g\in G(H)$, the linear space of $(h,g)${\it -primitive elements}
 is
$$
\mathcal{P}_{h,g}(H) :=\{x\in H\mid\com(x)= x\ot h + g\ot x\}.
$$
In case $g=1=h$, the linear space $\mathcal{P}(H) = \mathcal{P}_{1,1}(H)$ is called
the set of \textit{primitive elements}.

Let $M$ be
a left $H$-comodule via $\delta(m)=m_{(-1)}\ot m_{(0)} \in H\ot M$
for all $m\in M$. 
The space of {\it left
coinvariants} is given by $ ^{\co \delta}M = \{x\in
M\mid\delta(x)=1\ot x\}$. In particular, for a Hopf algebra map $\pi:H\rightarrow L$, 
it follows that
$H$ is a left $L$-comodule via $(\pi\ot\id)\com$ and
$$^{\co \pi}H:=\ ^{\co\ (\pi\ot\id)\com}H
=  \{h\in
H\mid (\pi\ot \id )\com(h)=1\ot h\}.$$ 
Right coinvariants, written
${H^{\co \pi}  }$, are defined analogously.

\subsection{Yetter-Drinfeld modules and
Nichols algebras}\label{subsec:ydCG-nichols-alg}

Let $H$ be a Hopf algebra. A \textit{left Yetter-Drinfeld module} $M$ over $H$
is a left $H$-module $(M,\cdot)$ and a left $H$-comodule 
$(M,\delta)$ satisfying 
$$ \delta(h\cdot m) = h_{(1)}m_{(-1)}\Ss(h_{(3)})\ot h_{(2)}\cdot m_{(0)}
\qquad \forall\ m\in M, h\in H. $$
We denote by $\ydh$ the category of left Yetter-Drinfeld modules over $H$. 
It is a braided monoidal category: for $M, N \in \ydh$, 
the braiding $c_{M,N}:M\ot N \to N\ot M$ is given by 
\begin{equation}\label{eq:braidingyd}
c_{M,N}(m\ot n) = m_{(-1)}\cdot n \ot m_{(0)} \qquad \forall\  
m\in M,n\in N. 
\end{equation}

\begin{defi}\cite[Definition 2.1]{AS} \label{defNicholsalgebra}
Let $H$ be a Hopf algebra and $V \in \ydh$. A braided $\N$-graded 
Hopf algebra $R = \bigoplus_{n\geq 0} R(n) \in \ydh$  is called 
the \textit{Nichols algebra} of $V$ if 
\begin{enumerate}
 \item[$(i)$] $R(0)\simeq \Bbbk$, $R(1)\simeq V$;
 \item[$(ii)$] $R(1) = \mathcal{P}(R)$;
 \item[$(iii)$] $R$ is generated as an algebra by $R(1)$.
\end{enumerate}
In this case, $R$ is denoted by $\toba(V) = \bigoplus_{n\geq 0} \toba^{n}(V) $.    
\end{defi}

For any $V \in \ydh$ there is a unique up to isomorphism Nichols algebra $\toba(V)$ associated with it. 
It is the quotient of the tensor algebra $T (V)$ by the largest 
homogeneous two-sided ideal $I$ satisfying:
\begin{itemize}
\item $I$ is generated by 
homogeneous elements of degree $\geq 2$, and
\item $\com(I) \subseteq I\ot T(V) + T(V)\ot I$, \textit{i.e.}, it is also a coideal.
\end{itemize}
See \cite[Section 2.1]{AS} for details.

\begin{obs}\label{def:nicholsbraiding}
Let $ c $ be the braiding associated to $ V\in\ydh $ and assume that there is $W\subseteq V$ a subspace such that $c(W\ot W) \subseteq W\ot W$. Then, one may identify
$\toba(W)$ with a subalgebra of $\toba(V)$; perhaps belonging to different braided monoidal categories.
In particular, $\toba(V)$ is infinite-dimensional whenever $\toba(W)$ is infinite-dimensional. 
This occurs for example, when 
$V$ contains a non-zero element $v$ such that $c(v\ot v) = v\ot v$.
\end{obs}

\subsection{Bosonization and Hopf algebras with a projection}\label{subsec:bosonization}
Let $H$ be a Hopf algebra and $B$ a braided Hopf algebra in $\ydh$.
The procedure to obtain a
usual Hopf algebra from $B$ and $H$ is called
the Majid-Radford biproduct or \emph{bosonization}, and it is usually
denoted by $B \#H$. As a vector space, $B \# H = B\otimes H$, and the
multiplication and comultiplication are given by the smash-product
and smash-coproduct,
respectively. 
Explicitly, for all $b, c \in B$ and $g,h \in H$, we have
\begin{align*}
(b \# g)(c \#h) & = b(g_{(1)}\cdot c)\# g_{(2)}h,\\
\com(b \# g) & =b^{(1)} \# (b^{(2)})_{(-1)}g_{(1)} \ot
(b^{(2)})_{(0)}\# g_{(2)},
\end{align*}
where $\com_{B}(b) = b^{(1)}\ot b^{(2)}$ denotes the comultiplication in $B\in \ydh$.
We identify $b=b\# 1$ and
$h=1\# h$; in particular we have $bh=b\# h$ and $hb=h_{(1)}\cdot b\# h_{(2)}$.
Clearly, the map $\iota: H \to B\#H$ given by $\iota(h) = 1\#h$ 
is an injective Hopf algebra map, and the map $\pi: B\#H \to H$
given by $\pi(b\#h) = \eps_{B}(b)h$ 
is a surjective Hopf algebra map such that $\pi \circ \iota = \id_{H} $.
Moreover, it holds that $B = (B\#H)^{\co \pi}$.

Conversely, let $A$ be a Hopf algebra with bijective antipode. Suppose that there are
Hopf algebra morphisms $\pi: A\to H$ and $\iota: H\to A$ 
such that $\pi\circ\iota =\id_{H}$.
Then $B=A^{\co\pi}$ is a braided Hopf algebra in $\ydh$ and $A\simeq B\# H$
as Hopf algebras.

\subsection{The Drinfeld double}\label{subsec:drinfeld-double}
We briefly describe the structure of the Drinfeld double of a finite-dimensional Hopf algebra.

Let $ H $ be a finite-dimensional Hopf algebra. 
Consider $H$ acting on $H^{*}$ and $H^{*}$ acting on $H$ via respectively
$$ h\twoheadrightarrow f = \langle f_{(3)} \Ss^{-1}(f_{(1)}), h
\rangle f_{(2)} \text{ and } h\twoheadleftarrow 
f = \langle f, \Ss^{-1}(h_{(3)})h_{(1)}\rangle h_{(2)}, \ \forall\ h\in H,\ f\in H^{*}. $$
The Drinfeld double of $H$ is the Hopf algebra $ D(H)$, 
where $ D(H) = H^{*}\otimes H$, 
as vector spaces. The product and the unit are given by
\[ (f\bowtie h)(g\bowtie k) = f(h_{(1)} \twoheadrightarrow g_{(2)})
\bowtie (h_{(2)}\twoheadleftarrow g_{(1)})k, \quad\text{and}\quad 1_{D(H)} =
\varepsilon\bowtie 1.
\]
The coproduct, counit and antipode do not play an important role in this paper. 
They can be found for instance in \cite[Definition 10.3.5]{Mo}.

The following result will be central in Section \ref{sec:categoryyd}, since we will
study the simple and indecomposable left $D(\Kk^{\cop})$-modules first and then 
translate the information to $\ydk$.

\begin{prop}\cite[Proposition 10.6.16]{Mo}
Let $ H $ be a finite-dimensional Hopf algebra. The category 
$ \ydh $ of left Yetter-Drinfeld modules over $H$ can be identified with the category $ _{D(H^{\cop})}\cal{M} $ of left 
modules over the Drinfeld double $ D(H^{\cop}) $.\qed
\end{prop}

\section{The Hopf algebra $\Kk$ and its Drinfeld double $D(\Kk^{\cop})$ }\label{sec:qsubgroupK}
%
All 
pointed nonsemisimple Hopf algebras of dimension $8$ were determined
by \c Stefan \cite{stefan}. Except for one case (up to isomorphism), these pointed Hopf algebras
have pointed duals.
The exception is given by 
$$
 {\mathcal A}''_4:=\Bbbk\langle g,x\mid g^4-1=x^2-g^2+1=gx+xg=0\rangle,
$$
with 
$\com(g)=g\ot g$ and $\com(x)=x\ot g+1\ot x$.
Moreover, it holds that $\mathcal{K} $, presented by \eqref{eq:defining_relations_K} and \eqref{eq:defining_coproduct_antipode_K}, is isomorphic to $ ({\mathcal A}''_{4})^{*}$, see \cite{GV}.
Up to isomorphism, $ \mathcal{K} $
is the only Hopf algebra of dimension $8$ which
is neither semisimple nor pointed nor has the Chevalley property.
The next proposition gives us a presentation
of $\Kk$ and some useful relations that will be used in the sequel.
The proof follows from
\cite[Lemma 3.3]{GV}.

Throughout the paper, we fix  a primitive $4$-th root of unity $\xi$.  

\begin{prop}\label{prop:desc-K} 
{\color{white} a}

\begin{enumerate}
  \item[$(i)$] $\mathcal{K}$ is generated as an algebra by the elements $a,b,c,d$ satisfying $ \eqref{eq:defining_relations_K} $.
\item[$(ii)$] The set $ \{1, a, b, c, d, a^2, ab, ac\}$ is a linear basis of $\Kk$.
\item[$(iii)$] The coalgebra structure and the antipode are determined by \eqref{eq:defining_coproduct_antipode_K}. In particular,
\begin{align}\label{eq:rest_of_coproduct_and_antipode_in_K}
 \Delta(ab) &= ab\otimes 1 + a^2\otimes ab, & \Delta(ac) &= ac\otimes a^2 + 1\otimes ac, & \Delta(a^2) &= a^2\otimes a^2, \\
\notag\Ss(a^2) &= a^2,& \Ss(ab) &= -ac,& \Ss(ac) &= ab.  
\end{align}
\item[$(iv)$]  The multiplication table of $\Kk$ is 
\begin{center}
\begin{tabular}{|c||c|c|c|c|c|c|c|c|}
\hline
& $ 1 $ & $ a $ & $ b $ & $ c $ & $ d $ & $ a^2 $ & $ ab $ & $ ac $ \\
\hline
\hline
$ 1 $ & $ 1 $ & $ a $ & $ b $ & $ c $ & $ d $ & $ a^2 $ & $ ab $ & $ ac $ \\
\hline
$ a $ & $ a $ & $ a^2 $ & $ ab $ & $ ac $ & $ 1 $ & $ d $ & $ c $ & $ b $ \\
\hline
$ b $ & $ b $ & $ -\xi ab $ & $ 0 $ & $ 0 $ & $ \xi ac $ & $ -c $ & $ 0 $ & $ 0 $ \\
\hline
$ c $ & $ c $ & $ -\xi ac $ & $ 0 $ & $ 0 $ & $ \xi ab $ & $ -b $ & $ 0 $ & $ 0 $ \\
\hline
$ d $ & $ d $ & $ 1 $ & $ ac $ & $ ab $ & $ a^2 $ & $ a $ & $ b $ & $ c $ \\
\hline
$ a^2 $ & $ a^2 $ & $ d $ & $ c $ & $ b $ & $ a $ & $ 1 $ & $ ac $ & $ ab $ \\
\hline
$ ab $ & $ ab $ & $ -\xi c $ & $ 0 $ & $ 0 $ & $ \xi b $ & $ -ac $ & $ 0 $ & $ 0 $ \\
\hline
$ ac $ & $ ac $ & $ -\xi b $ & $ 0 $ & $ 0 $ & $ \xi c $ & $ -ab $ & $ 0 $ & $ 0 $ \\
\hline
\end{tabular}
\end{center}
\item[$(v)$] $\mathcal{K} \simeq H_{4}\oplus \mathcal M^{\ast}(2,\Bbbk)$
as coalgebras, where $H_{4}$ is the Sweedler's Hopf algebra and 
$\mathcal{M}^{\ast}(2,\Bbbk)$ is a comatrix coalgebra of dimension $4$.
 \end{enumerate}\qed
\end{prop}

\begin{rmks}\label{rmk:strucdual}
	{\color{white} a}

\noindent $(a)$ Denote by $ \{1^*, a^*, b^*, c^*, d^*, (a^2)^*, (ab)^*, (ac)^*\} $ the 
basis of $ \Kk^* $ dual to  $ \{1, a, b, c, d, a^2, ab, $  $ ac\} $.
Using the multiplication table in Proposition \ref{prop:desc-K} (iv), it follows that
\begin{align*}
\Delta(1^*) &= 1^*\otimes 1^* + a^*\otimes d^* + d^*\otimes a^* + (a^2)^*\otimes (a^2)^* ,\\
\Delta(a^*) &= 1^*\otimes a^* + a^*\otimes 1^* + (a^2)^*\otimes d^* + d^*\otimes (a^2)^* ,\\
\Delta(d^*) &= 1^*\otimes d^* + d^*\otimes 1^* + (a^2)^*\otimes a^* + a^*\otimes (a^2)^* ,\\
\Delta((a^2)^*) &= 1^*\otimes (a^2)^* + (a^2)^*\otimes 1^* + a^*\otimes a^* + d^*\otimes d^*, \\
\Delta(b^*) &= 1^*\otimes b^* + b^*\otimes 1^* + a^*\otimes (ac)^* -\xi (ac)^*\otimes a^*+ \\
&+ (a^2)^*\otimes c^* - c^*\otimes (a^2)^* +\xi (ab)^*\otimes d^* + d^*\otimes (ab)^* ,\\
\Delta(c^*) &= 1^*\otimes c^* + c^*\otimes 1^* -\xi (ab)^*\otimes a^* + a^*\otimes (ab)^* +\\
&+ (a^2)^*\otimes b^* - b^*\otimes (a^2)^* +\xi (ac)^*\otimes d^* + d^*\otimes (ac)^* ,\\
\Delta((ab)^*) &= 1^*\otimes (ab)^* + (ab)^*\otimes 1^* -\xi b^*\otimes a^* + a^*\otimes b^* +\\
&+ d^*\otimes c^* +\xi c^*\otimes d^* - (ac)^*\otimes (a^2)^* + (a^2)^*\otimes (ac)^* ,\\
\Delta((ac)^*) &= 1^*\otimes (ac)^* + (ac)^*\otimes 1^* -\xi c^*\otimes a^* + a^*\otimes c^* +\\
&+ d^*\otimes b^* +\xi b^*\otimes d^* - (ab)^*\otimes (a^2)^* + (a^2)^*\otimes (ab)^*.
\end{align*}

\noindent
$(b)$ Let $ \alpha\in G(\Kk^*) = \Alg({\cal{K}}, \Bbbk)$. 
The relations \eqref{eq:defining_relations_K} implies that
$\alpha(a) $ is a $4$-th root of unity,  $ \alpha(b) = \alpha(c)=0 $ and $ \alpha(d) = \alpha(a)^{-1}$. 
Thus $G(\Kk^*) $ consists of the elements
\begin{align*}
\alpha_{j} = 1^* +  \xi^{-j}a^* + \xi^{j} d^* + (-1)^{j} (a^2)^*,\qquad  j=0,1,2,3.
\end{align*}
Note that $ \alpha_0 = \varepsilon $ and $\alpha_{1}^{j}=\alpha_{j}$. In particular,
$G(\Kk^*)\simeq \Z / 4\Z$ and 
$ \alpha_1, \alpha_{3} $ 
are generators. 

\noindent
$(c)$
The multiplication table of $\Kk^*$ is 
\begin{center}
\begin{tabular}{|c||c|c|c|c|c|c|c|c|}
\hline
& $ 1^{*} $ & $ a^{*} $ & $ b^{*} $ & $ c^{*} $ & $ d^{*} $ & $ (a^2)^{*} $ & $ (ab)^{*} $ & $ (ac)^{*} $ \\
\hline
\hline
$ 1^{*} $ & $ 1^{*} $ & $ 0 $ & $ 0 $ & $ 0 $ & $ 0 $ & $ 0 $ & $ 0 $ & $ (ac)^{*} $ \\
\hline
$ a^{*} $ & $ 0 $ & $ a^{*} $ & $ b^{*} $ & $ 0 $ & $ 0 $ & $ 0 $ & $ 0 $ & $ 0 $ \\
\hline
$ b^{*} $ & $ 0 $ & $ 0 $ & $ 0 $ & $ a^{*} $ & $ b^{*} $ & $ 0 $ & $ 0 $ & $ 0 $ \\
\hline
$ c^{*} $ & $ 0 $ & $ c^{*} $ & $ d^{*} $ & $ 0 $ & $ 0 $ & $ 0 $ & $ 0 $ & $ 0 $ \\
\hline
$ d^{*} $ & $ 0 $ & $ 0 $ & $ 0 $ & $ c^{*} $ & $ d^{*} $ & $ 0 $ & $ 0 $ & $ 0 $ \\
\hline
$ (a^2)^{*} $ & $ 0 $ & $ 0 $ & $ 0 $ & $ 0 $ & $ 0 $ & $ (a^2)^{*} $ & $ (ab)^{*} $ & $ 0 $ \\
\hline
$ (ab)^{*} $ & $ (ab)^{*} $ & $ 0 $ & $ 0 $ & $ 0 $ & $ 0 $ & $ 0 $ & $ 0 $ & $ 0 $ \\
\hline
$ (ac)^{*} $ & $ 0 $ & $ 0 $ & $ 0 $ & $ 0 $ & $ 0 $ & $ (ac)^{*} $ & $ 0 $ & $ 0 $ \\
\hline
\end{tabular}
\end{center}
\end{rmks}

In order to compute the Drinfeld double $D(\Kk^{\cop})$, we need to describe 
the isomorphism $\Kk^{*} \simeq \Aa$ explicitly.

\begin{lema}\label{lem:isomAK}
The algebra map $\varphi: \Aa \to \Kk^{*}$ given by 
$$\varphi(g) = \alpha_{1} \quad\text{ and }\quad
\varphi(x) =\sqrt{2}\xi(b^* + c^* + (ab)^* + (ac)^*),
$$
is a Hopf algebra isomorphism.
\end{lema}

\pf  
A direct computation shows that $\varphi$ is a coalgebra map.
Hence, the image of $\varphi$ is a Hopf subalgebra of $\Kk^{*}$ of dimension greater than $4$, 
because it contains
the group algebra $\Bbbk G(\Kk^{*})$ and the image of the skew-primitive element $x$. 
By the Nichols-Zoeller theorem it follows that $\varphi$ is surjective and whence 
an isomorphism.
\epf

\begin{obs}
Consider the basis 
$ \{g^j, xg^j\}_{ 0\leq j \leq 3} $ of $\mathcal{A}_{4}''$. 
By Remark \ref{rmk:strucdual} (c) and Lemma 
\ref{lem:isomAK}, 
it follows that 
\begin{align*}
\varphi(g^j) &= \alpha_{j}  &\text{ for all } & &  0\leq j\leq 3, \\
\varphi(xg^j) &=\sqrt{2}\xi(\xi^j b^* +
\xi^{-j} c^* + (ab)^* +(-1)^j(ac)^*) & \text{ for all } & &0\leq j\leq 3.
\end{align*} 
\end{obs}

\subsection{Description of $D(\Kk^{\cop})$}
In this subsection we describe the Drinfeld double $D(\Kk^{\cop})$. To make the
notation lighter, from now on we write $D = D(\Kk^{\cop}) $.

\begin{prop}\label{prop:Dgen-rel}
$D$ is the $\Bbbk$-algebra generated by the elements  $a, b, c, d, x, g$ such that 
$a,b,c,d$
satisfy the relations of $\Kk^{\cop}$;
$x,g$ satisfy the relations of $(\Aa)^{\op \cop} $; and all together they
satisfy the following relations:
\begin{tabbing}
\hspace{1cm}\=\hspace{7cm}\=\kill
\>  $ ax +\xi xa = \sqrt{2}\xi(b + gc) $, \>  $ bx -\xi xb = \sqrt{2}\xi(a - gd) $, \\
\> \ \ \ \ \ \ \ \ \ \ $ ag = ga $, \> \ \ \ \ \ \ \ \ \ \ $ bg = -gb $, \\
\> \ \ \ \ \ \ \ \ \ \ $ cg = - gc $, \> \ \ \ \ \ \ \ \ \ \ $ dg = gd $,\\
\> $ cx + \xi xc = \sqrt{2}\xi(d - ga) $, \>  $ dx - \xi xd = \sqrt{2}\xi(c + gb) $. 
\end{tabbing}
\end{prop}

\pf 
Since $ (f\bowtie 1)(g\bowtie k) = fg\bowtie k $ and 
$ (f\bowtie h)(1\bowtie k) = f\bowtie hk$ for all $f, g\in 
(\Aa)^{\op \cop}$ and $h, k \in \Kk^{\cop} $, 
it is enough to describe the relations derived from products of the form
$(1_{\Aa}\bowtie h)(y\bowtie 1_{\Kk})$,
 where $ h\in \Kk^{\cop}$ and $ y\in (\Aa)^{\op \cop}$ are algebra generators. 

Assume first that $ y = g$. Since
$h\twoheadrightarrow g = \langle \Ss_{\cop}^{-1}(g) g, h\rangle g =
\langle \Ss(g) g, h\rangle g = \langle 1, h\rangle g = \varepsilon(h)g, 
$
for all  $h\in \Kk^{\cop}$, it follows that 
\begin{align*}
& (1_{\Aa}\bowtie 
\left\{
\begin{array}{c}
a \\
b \\
c \\
d
\end{array}
\right\})(g\bowtie 1_{\Kk}) = \left\{
\begin{array}{c}
{a_{(1)}}_{\cop}\twoheadrightarrow g_{(2)}\bowtie {a_{(2)}}_{\cop}\twoheadleftarrow g_{(1)} \\
{b_{(1)}}_{\cop}\twoheadrightarrow g_{(2)}\bowtie {b_{(2)}}_{\cop}\twoheadleftarrow g_{(1)} \\
{c_{(1)}}_{\cop}\twoheadrightarrow g_{(2)}\bowtie {c_{(2)}}_{\cop}\twoheadleftarrow g_{(1)} \\
{d_{(1)}}_{\cop}\twoheadrightarrow g_{(2)}\bowtie {d_{(2)}}_{\cop}\twoheadleftarrow g_{(1)}
\end{array}
\right\} \\
&= \left\{
\begin{array}{c}
a\twoheadrightarrow g \bowtie a\twoheadleftarrow g +
c\twoheadrightarrow g \bowtie  b\twoheadleftarrow g \\
b\twoheadrightarrow g \bowtie a\twoheadleftarrow g + 
d\twoheadrightarrow g \bowtie 
b\twoheadleftarrow g \\
a\twoheadrightarrow g \bowtie c\twoheadleftarrow g + 
c\twoheadrightarrow g \bowtie d\twoheadleftarrow g \\
b\twoheadrightarrow g \bowtie c\twoheadleftarrow g +
d\twoheadrightarrow g \bowtie d\twoheadleftarrow g
\end{array}
\right\} = g\bowtie \left\{
\begin{array}{c}
a\twoheadleftarrow g \\
b\twoheadleftarrow g \\
c\twoheadleftarrow g \\
d\twoheadleftarrow g
\end{array}
\right\}\\
&= g\bowtie \left\{
\begin{array}{c}
\langle g, \Ss_{\cop}^{-1}({a_{(3)}}_{\cop}){a_{(1)}}_{\cop}\rangle {a_{(2)}}_{\cop} \\
\langle g, \Ss_{\cop}^{-1}({b_{(3)}}_{\cop}){b_{(1)}}_{\cop}\rangle {b_{(2)}}_{\cop} \\
\langle g, \Ss_{\cop}^{-1}({c_{(3)}}_{\cop}){c_{(1)}}_{\cop}\rangle {c_{(2)}}_{\cop} \\
\langle g, \Ss_{\cop}^{-1}({d_{(3)}}_{\cop}){d_{(1)}}_{\cop}\rangle {d_{(2)}}_{\cop}
\end{array}
\right\} = g\bowtie \left\{
\begin{array}{c}
\langle g, \Ss(a_{(1)})a_{(3)}\rangle a_{(2)} \\
\langle g, \Ss(b_{(1)})b_{(3)}\rangle b_{(2)} \\
\langle g, \Ss(c_{(1)})c_{(3)}\rangle c_{(2)} \\
\langle g, \Ss(d_{(1)})d_{(3)}\rangle d_{(2)}
\end{array}
\right\} \\
&= g\bowtie \left\{
\begin{array}{c}
\langle g, \Ss(a)a\rangle a + \langle g, \Ss(b)a\rangle c +\langle g, \Ss(a)c\rangle b +\langle g, \Ss(b)c\rangle d \\
\langle g, \Ss(a)b\rangle a + \langle g, \Ss(b)b\rangle c +\langle g, \Ss(a)d\rangle b +\langle g, \Ss(b)d\rangle d \\
\langle g, \Ss(c)a\rangle a + \langle g, \Ss(d)a\rangle c +\langle g, \Ss(c)c\rangle b +\langle g, \Ss(d)c\rangle d \\
\langle g, \Ss(c)b\rangle a + \langle g, \Ss(d)b\rangle c +\langle g, \Ss(c)d\rangle b +\langle g, \Ss(d)d\rangle d \\
\end{array}
\right\} \\
&= g\bowtie \left\{
\begin{array}{c}
\langle g, da\rangle a + \langle g, \xi ba\rangle c +\langle g, dc\rangle b +\langle g, \xi bc\rangle d \\
\langle g, db\rangle a + \langle g, \xi bb\rangle c +\langle g, dd\rangle b +\langle g, \xi bd\rangle d \\
\langle g, -\xi ca \rangle a + \langle g, aa\rangle c +\langle g, -\xi cc\rangle b +\langle g, ac\rangle d \\
\langle g, -\xi cb\rangle a + \langle g, ab\rangle c +\langle g, -\xi cd\rangle b +\langle g, ad\rangle d \\
\end{array}
\right\} \\
&= g\bowtie \left\{
\begin{array}{c}
\langle g, 1\rangle a + \langle g, ab\rangle c +\langle g, ab\rangle b +\langle g, 0\rangle d \\
\langle g, ac\rangle a + \langle g, 0\rangle c +\langle g, a^2\rangle b +\langle g, -ac\rangle d \\
\langle g, -ac\rangle a + \langle g, a^2\rangle c +\langle g, 0\rangle b +\langle g, ac\rangle d \\
\langle g, 0\rangle a + \langle g, ab\rangle c +\langle g, ab\rangle b +\langle g, 1\rangle d \\
\end{array}
\right\} = g\bowtie \left\{
\begin{array}{c}
a \\
-b \\
-c \\
d \\
\end{array}
\right\}.
\end{align*}
From this equalities it follows that $ ag = ga, bg = -gb, cg = -gc $ and $ dg = gd$.

Suppose now that 
$ y = x $. 
Using the computations above, we have that
\begin{align*}
\left\{
\begin{array}{c}
a \\
b \\
c \\
d
\end{array}\right\} \twoheadleftarrow x  = \left\{
\begin{array}{c}
\langle x, 1\rangle a + \langle x, ab\rangle c +\langle x, ab\rangle b +\langle x, 0\rangle d \\
\langle x, ac\rangle a + \langle x, 0\rangle c +\langle x, a^2\rangle b +\langle x, -ac\rangle d \\
\langle x, -ac\rangle a + \langle x, a^2\rangle c +\langle x, 0\rangle b +\langle x, ac\rangle d \\
\langle x, 0\rangle a + \langle x, ab\rangle c +\langle x, ab\rangle b +\langle x, 1\rangle d \\
\end{array}
\right\} = \sqrt{2}\xi \left\{
\begin{array}{c}
b+c \\
a-d \\
d-a \\
b+c \\
\end{array}
\right\},
\end{align*}
and
\begin{align*}
\left\{
\begin{array}{c}
a \\
b \\
c \\
d
\end{array}
\right\}\twoheadrightarrow x &= \langle\Ss(x_{(1)})x_{(3)}, \left\{
\begin{array}{c}
a \\
b \\
c \\
d \\
\end{array}
\right\} \rangle x_{(2)} \\
&= 
\langle\Ss(x)g, \left\{
\begin{array}{c}
a \\
b \\
c \\
d \\
\end{array}
\right\} \rangle g + \langle\Ss(1)g, \left\{
\begin{array}{c}
a \\
b \\
c \\
d \\
\end{array}
\right\} \rangle x + \langle\Ss(1)x, \left\{
\begin{array}{c}
a \\
b \\
c \\
d \\
\end{array}
\right\} \rangle 1\\
&= \langle -xg^3 g, \left\{
\begin{array}{c}
a \\
b \\
c \\
d \\
\end{array}
\right\} \rangle g + \langle g, \left\{
\begin{array}{c}
a \\
b \\
c \\
d \\
\end{array}
\right\} \rangle x + \langle x, \left\{
\begin{array}{c}
a \\
b \\
c \\
d \\
\end{array}
\right\} \rangle 1 \\
&= \langle x, \left\{
\begin{array}{c}
a \\
b \\
c \\
d \\
\end{array}
\right\} \rangle (1 - g) + \langle g, \left\{
\begin{array}{c}
a \\
b \\
c \\
d \\
\end{array}
\right\} \rangle x = \left\{
\begin{array}{c}
\xi x \\
\sqrt{2}\xi (1 - g) \\
\sqrt{2}\xi (1 - g) \\
-\xi x 
\end{array}
\right\}.
\end{align*}
Hence, it follows that 
\begin{align*}
& (1_{\Aa}\bowtie 
\left\{
\begin{array}{c}
a \\
b \\
c \\
d
\end{array}
\right\})(x\bowtie 1_{\Kk}) = \left\{
\begin{array}{c}
a_{(2)}\twoheadrightarrow g \bowtie a_{(1)}\twoheadleftarrow x +
a_{(2)}\twoheadrightarrow x \bowtie a_{(1)}\twoheadleftarrow 1 \\
b_{(2)}\twoheadrightarrow g \bowtie b_{(1)}\twoheadleftarrow x +
b_{(2)}\twoheadrightarrow x \bowtie b_{(1)}\twoheadleftarrow 1 \\
c_{(2)}\twoheadrightarrow g \bowtie c_{(1)}\twoheadleftarrow x +
c_{(2)}\twoheadrightarrow x \bowtie c_{(1)}\twoheadleftarrow 1 \\
d_{(2)}\twoheadrightarrow g\bowtie d_{(1)}\twoheadleftarrow x +
d_{(2)}\twoheadrightarrow x \bowtie d_{(1)}\twoheadleftarrow 1
\end{array}
\right\} 
\\
&= \left\{
\begin{array}{c}
\varepsilon(a_{(2)}) g \bowtie a_{(1)}\twoheadleftarrow x +
a_{(2)}\twoheadrightarrow x \bowtie a_{(1)} \\	
\varepsilon(b_{(2)}) g \bowtie b_{(1)}\twoheadleftarrow x +
b_{(2)}\twoheadrightarrow x \bowtie b_{(1)} \\
\varepsilon(c_{(2)}) g \bowtie c_{(1)}\twoheadleftarrow x + 
c_{(2)}\twoheadrightarrow x \bowtie c_{(1)} \\
\varepsilon(d_{(2)}) g \bowtie d_{(1)}\twoheadleftarrow x + 
d_{(2)}\twoheadrightarrow x \bowtie d_{(1)}
\end{array}
\right\} \\
&= \left\{
\begin{array}{c}
g \bowtie  a\twoheadleftarrow x +
a\twoheadrightarrow x \bowtie a +
c\twoheadrightarrow x  \bowtie b\\
g \bowtie  b\twoheadleftarrow x +
b\twoheadrightarrow x \bowtie a +
d\twoheadrightarrow x  \bowtie b\\
g \bowtie  c\twoheadleftarrow x +
a\twoheadrightarrow x \bowtie c+
c\twoheadrightarrow x  \bowtie d\\
g \bowtie  d\twoheadleftarrow x +
b\twoheadrightarrow x \bowtie c +
d\twoheadrightarrow x  \bowtie d
\end{array} \right\} \\
&= \left\{
\begin{array}{c}
\sqrt{2}\xi g \bowtie (b+c) +
\xi x \bowtie a +
\sqrt{2}\xi (1 - g) \bowtie b\\
\sqrt{2}\xi g \bowtie (a-d) +
\sqrt{2}\xi (1 - g) \bowtie a 
-\xi x \bowtie b\\
\sqrt{2}\xi g \bowtie (d-a) +
\xi x \bowtie c +
\sqrt{2}\xi (1 - g) \bowtie d\\
\sqrt{2}\xi g \bowtie (b+c) +
\sqrt{2}\xi (1 - g) \bowtie c 
-\xi x \bowtie d
\end{array}
\right\}\\
&= \left\{
\begin{array}{c}
-\xi x \bowtie a +
\sqrt{2}\xi g\bowtie c
+ \sqrt{2}\xi \bowtie b\\
\xi x \bowtie b 
- \sqrt{2}\xi g\bowtie d
+ \sqrt{2}\xi \bowtie a\\
-\xi x \bowtie c -
 \sqrt{2}\xi g\bowtie a
+ \sqrt{2}\xi \bowtie d\\
\xi x \bowtie d +
 \sqrt{2}\xi g\bowtie b
+ \sqrt{2}\xi \bowtie c
\end{array}
\right\},
\end{align*}
which gives us the other four relations of $ D$.
\epf

\subsection{Simple left $D$-modules}
We begin by describing the $1$-dimensional  $ D$-modules.
Given a character $\chi$ on $D$, we denote by $\Bbbk_{\chi}$ the module associated with it.

\begin{lema} \label{modulos dim 1}
There are four non-isomorphic $1$-dimensional left $ D$-modules given by
the characters $\chi^{j},\ 0\leq j\leq 3 $, where
$$\chi^{j}(a) = \xi^{j}, \ \chi^{j}(b) = 0, \ \chi^{j}(c) = 0, \ \chi^{j}(d) = 
\xi^{-j}, \ \chi^{j}(x) = 0, \ \chi^{j}(g) = (-1)^{j}.
$$
Moreover, any $1$-dimensional $ D$-module is isomorphic to $\Bbbk_{\chi^{j}}$
for some $0\leq j\leq 3$.
\end{lema}

\pf Straightforward.
\epf

We describe next the simple $ D$-modules of dimension two.
For this, consider the set  
$$\Lambda = \{(i,j) \in \Z_{4}\times \Z_{4}|\  
2i\neq j\}.$$
Clearly, $|\Lambda| = 12$.

\begin{lema} \label{modulos dim 2}
For any pair $(i,j) \in \Lambda$,
there exists a simple $2$-dimensional left $ D$-module 
$V_{i,j}$. 
The action on a fixed basis  is given by 
\[
\rho_{i,j}(a)= \left(
\begin{array}{cc}
\xi^{i} & 0 \\
0 & \xi^{i+3} \\
\end{array}
\right), \ 
\rho_{i,j}(b)= \left(
\begin{array}{cc}
0 & (-1)^{i} \\
0 & 0 \\
\end{array}
\right), \ 
\rho_{i,j}(c)= \left(
\begin{array}{cc}
0 & 1 \\
0 & 0 \\
\end{array}
\right),
\]\[
\rho_{i,j}(d)= \left(
\begin{array}{cc}
\xi^{-i} & 0 \\
0 & \xi^{-i+1} \\
\end{array}
\right), \quad
\rho_{i,j}(g)= \left(
\begin{array}{cc}
\xi^{j} & 0 \\
0 & \xi^{j+2} \\
\end{array}
\right),
 \]\[
\rho_{i,j}(x)= \left(
\begin{array}{cc}
0 & \dfrac{\sqrt{2}}{2}\xi(\xi^{i} + \xi^{3i+j}) \\
\sqrt{2}\xi(\xi^{3i} - \xi^{i+j}) & 0 \\
\end{array}
\right),
\]
Moreover, any simple $2$-dimensional $ D$-module is 
isomorphic to $V_{i,j}$ for some $(i,j)\in \Lambda$, and 
$V_{i,j}\simeq V_{k,\ell}$ if and only if
$(i,j)  = (k,\ell)$.
\end{lema}

\pf Let $\rho: D \to \End(V)$ be a $2$-dimensional simple representation  and assume that the 
associated matrices of the generators of $D$ on a fixed basis of $V$ are given by:
\begin{align*}
\rho(a)= \left(
\begin{array}{cc}
a_{11} & a_{12} \\
a_{21} & a_{22} \\
\end{array}
\right), \ \rho(b)= \left(
\begin{array}{cc}
b_{11} & b_{12} \\
b_{21} & b_{22} \\
\end{array}
\right), \ \rho(c)= \left(
\begin{array}{cc}
c_{11} & c_{12} \\
c_{21} & c_{22} \\
\end{array}
\right), \\
\rho(d)= \left(
\begin{array}{cc}
d_{11} & d_{12} \\
d_{21} & d_{22} \\
\end{array}
\right), \ \rho(x)= \left(
\begin{array}{cc}
x_{11} & x_{12} \\
x_{21} & x_{22} \\
\end{array}
\right), \ \rho(g)= \left(
\begin{array}{cc}
g_{11} & g_{12} \\
g_{21} & g_{22} \\
\end{array}
\right).
\end{align*}
As $ a^4 = 1=g^{4}$ and $ga=ag$, $\rho(a) $ and  
$\rho(g)$ are simultaneously diagonalizable and, without loss of generality, we may assume that 
\begin{align*}
\rho(a)= \left(
\begin{array}{cc}
\lambda_1 & 0 \\
0 & \lambda_2 \\
\end{array}
\right),\quad 
\rho(d)= \left(
\begin{array}{cc}
\lambda_1^{-1} & 0 \\
0 & \lambda_2 ^{-1}\\
\end{array}
\right)
\text{ and } \rho(g)= \left(
\begin{array}{cc}
\lambda_3 & 0 \\
0 & \lambda_4 \\
\end{array}
\right),
\end{align*}
where $ \lambda_i^4 = 1$ for $1\leq i \leq 4$. 
From the relation $ ac = \xi ca$ we have that
\begin{align*}
\left(
\begin{array}{cc}
\lambda_1c_{11} & \lambda_1c_{12} \\
\lambda_2c_{21} & \lambda_2c_{22} \\
\end{array}
\right) = \xi \left(
\begin{array}{cc}
\lambda_1c_{11} & \lambda_2c_{12} \\
\lambda_1c_{21} & \lambda_2c_{22} \\
\end{array}
\right),
\end{align*}
which implies that $ c_{11} = c_{22} = 0$. 
Similarly, the relation $ gx = - xg $ implies $ x_{11} = x_{22} = 0 $. 
Since $ a^2c = b$, we must have that
\begin{align*}
\rho(b) = \left(
\begin{array}{cc}
0 & \lambda_1^2c_{12} \\
\lambda_2^2c_{21} & 0 \\
\end{array}
\right).
\end{align*}
Also note that from the relation $ c^2 = 0$, we get that
$c_{12}c_{21} = 0$.
Thus, by permuting the elements of the basis, we may assume that $ c_{21} = 0$.
Suppose $ c_{12} = 0. $ That is,
\[
\rho(b) = \rho(c) = \left(
\begin{array}{cc}
0 & 0 \\
0 & 0 \\
\end{array}
\right) \text{ and } 
\rho(x)= \left(
\begin{array}{cc}
0 & x_{12} \\
x_{21} & 0 \\
\end{array}
\right).\] 
Clearly, these modules are simple if and only if $ x_{12} \neq 0 $ and $ x_{21}\neq 0$. 
As $ ax = -\xi xa + \sqrt{2}\xi(b + gc), $ it follows that
$x_{12}(\lambda_1 + \xi\lambda_2) = 0 $ and 
$x_{21}(\lambda_2 + \xi\lambda_1) = 0$.
Since $ x_{12}x_{21}\neq 0$, we must have that $ \lambda_1 + \xi\lambda_2 = \lambda_2 + \xi\lambda_1 = 0$,
which implies that $ \lambda_1 = \lambda_2 = 0$, a contradiction.

Therefore, we must have that $ c_{12}\neq 0$.
Clearly, we may also assume that $ c_{12} = 1$.
From the equality $ ac = \xi ca$, we get that $ \lambda_2 = -\xi\lambda_1$.
Moreover, seeing that $ cg = -gc, $ we must have $ \lambda_4 = -\lambda_3 $.  
Now, the relation $ ax +\xi xa = \sqrt{2}\xi(b +gc)$ yields
\begin{align*}
\left(
\begin{array}{cc}
0 & \lambda_1 x_{12} \\
\lambda_2 x_{21} & 0 \\
\end{array}
\right) &= -\xi \left(
\begin{array}{cc}
0 & \lambda_2 x_{12} \\
\lambda_1 x_{21} & 0 \\
\end{array}
\right) + \sqrt{2}\xi \left(
\begin{array}{cc}
0 & \lambda_1^2 + \lambda_3 \\
0 & 0 \\
\end{array}
\right),
\end{align*}
which implies that $ x_{12} = \dfrac{\sqrt{2}}{2}\xi(\lambda_1 + \lambda_3\lambda_1^{-1}) $. 
This is the same information obtained from the relation
$ dx - \xi xd = \sqrt{2}\xi(c + gb) $. Analogously, $ cx + \xi xc =\sqrt{2}\xi(d - ga) $ yields
\begin{align*}
\left(
\begin{array}{cc}
x_{21} & 0 \\
0 & \xi x_{21} \\
\end{array}
\right) = \sqrt{2}\xi \left(
\begin{array}{cc}
\lambda_1^{-1} - \lambda_1\lambda_3 & 0 \\
0 & \xi(\lambda_1^{-1} - \lambda_1\lambda_3) \\
\end{array}
\right),
\end{align*}
which gives $ x_{21} = \sqrt{2}\xi(\lambda_1^{-1} - \lambda_1\lambda_3)$.  
Also, the relation $ bx -\xi xb = \sqrt{2}\xi(a-gd)$ yields no further condition on the coefficients. 
Considering $ g^2 = 1 + x^2 $, we must have that $ x_{12}x_{21} = \lambda_3^2 - 1 $. 
In fact,
$$
x_{12}x_{21} = \dfrac{\sqrt{2}}{2}\xi(\lambda_1 +
\lambda_3\lambda_1^{-1})\sqrt{2}\xi(\lambda_1^{-1} - 
\lambda_1\lambda_3) = -(1 -\lambda_1^2\lambda_3+\lambda_3\lambda_1^{-2} 
-\lambda_3^2) 
= \lambda_3^2 - 1.$$

From the discussion above, the matrices defining the action on the simple module $V$ 
are of the form
 \[
\rho(a)= \left(
\begin{array}{cc}
\lambda_1 & 0 \\
0 & -\xi\lambda_1 \\
\end{array}
\right), \ \rho(b)= \left(
\begin{array}{cc}
0 & \lambda_1^2 \\
0 & 0 \\
\end{array}
\right), \ \rho(c)= \left(
\begin{array}{cc}
0 & 1 \\
0 & 0 \\
\end{array}
\right), \
\rho(g)= \left(
\begin{array}{cc}
\lambda_3 & 0 \\
0 & -\lambda_3 \\
\end{array}
\right),
\]\[
\rho(d)= \left(
\begin{array}{cc}
\lambda_1^{-1} & 0 \\
0 & \xi\lambda_1^{-1} \\
\end{array}
\right), \ \rho(x)= \left(
\begin{array}{cc}
0 & \dfrac{\sqrt{2}}{2}\xi(\lambda_1 +\lambda_3\lambda_1^{-1}) \\
\sqrt{2}\xi(\lambda_1^{-1} - \lambda_1\lambda_3) & 0 \\
\end{array}
\right), \]
with $ \lambda_1^4 = 1 =\lambda_3^4$. Moreover,
a direct computation shows that $V$ is simple if and only if
$\lambda_3 \neq \lambda_1^2$. Set $\lambda_{1}=\xi^{i}$ and $\lambda_{3}= \xi^{j}$ for some
$i,j\in \Z_{4}$. Then $2i\neq j$ and consequently 
$(i,j) \in \Lambda$.

\smallbreak
Finally, we show that $V_{i,j}$ is isomorphic to 
$V_{k,\ell}$ if and only if $(i,j) =
(k,\ell)$.
Let $ T: V_{i,j} \rightarrow V_{k,\ell} $
be an isomorphism of $ D$-modules; \textit{i.e.},
$\rho_{k,\ell}(t) T = T\rho_{i,j}(t)$ for all $t\in  D$. Denote by $[T] = \left(
\begin{array}{cc}
t_{11} & t_{12} \\
t_{21} & t_{22} \\
\end{array}
\right) $ the matrix of $T$ with respect to the given basis. 
Using the action of $ c $, we must have that
$ t_{21} = 0 $ and $ t_{11} = t_{22}$, because
\[ \left(
\begin{array}{cc}
t_{21} & t_{22} \\
0 & 0 \\
\end{array}
\right) = \left(
\begin{array}{cc}
0 & 1 \\
0 & 0 \\
\end{array}
\right)\left(
\begin{array}{cc}
t_{11} & t_{12} \\
t_{21} & t_{22} \\
\end{array}
\right) = \left(
\begin{array}{cc}
t_{11} & t_{12} \\
t_{21} & t_{22} \\
\end{array}
\right)\left(
\begin{array}{cc}
0 & 1 \\
0 & 0 \\
\end{array}
\right) = \left(
\begin{array}{cc}
0 & t_{11} \\
0 & t_{21} \\
\end{array}
\right).
\]
 Moreover, acting by $ a$ we obtain that
\[ \left(
\begin{array}{cc}
\xi^{k}t_{11} & \xi^{k}t_{12} \\
0 & -\xi^{k+1}t_{11} \\
\end{array}
\right) = 
[\rho_{k,\ell}(a)] [T] = [T][\rho_{i,j}(a)]
= \left(
\begin{array}{cc}
\xi^{i}t_{11} & -\xi^{i+1}t_{12} \\
0 & -\xi^{i+1}t_{11} \\
\end{array}
\right),
\]
which implies $ (\xi^{k} - \xi^{i})t_{11} = 0 $ and 
$ (\xi^{k} + \xi^{i+1})t_{12} = 0$. 
Since $ T $ is an isomorphism, 
this implies that $\xi^{k} = \xi^{i}$, from which follows that
$t_{12} = 0$. Consequently, $[T] = t_{11}I$.
Finally, acting by $ g$ yields that $ \xi^{\ell} = \xi^{j}$ and the claim
follows.
\epf

\begin{obs}\label{rmk:dualsVij}
Let $V$ be a left $D$-module. Since $D$ is a Hopf algebra, $V^{*}$ inherits a left $D$-module 
structure by the formula $(h\cdot f)(v)= f(\Ss(h)\cdot v)$ for all $f\in V^{*}$, $v\in V$ and 
$h\in D$. A straightforward computation yields $V_{i,j}^{*} \simeq V_{-i+1, -j+2}$ for all 
$(i,j) \in \Lambda$, where the indices in the second term are considered modulo $4$.
\end{obs}

We end this subsection by describing all simple left $ D $-modules up to isomorphism. 

\begin{teor}\label{thm:simpleD-mod}
There are sixteen simple left  $ D $-modules pairwise non-isomorphic: four
$1$-dimensional ones, given by Lemma \ref{modulos dim 1}, and twelve 
$2$-dimensional ones, given by Lemma \ref{modulos dim 2}.
\end{teor}

\pf Assume that there is a
simple module of dimension $d>2$ and let $n$
be the amount of simple $d$-dimensional modules 
pairwise non-isomorphic. Since $D$ is non-semisimple, 
the projective covers of the $1$-dimensional modules
have dimension at least $2$. Thus, by
 Lemmata \ref{modulos dim 1} and \ref{modulos dim 2}, 
it follows that
$$4.2 + 12.2^2 + nd^{2}= 56 + nd^{2}\leq 
\dim{D} = 64.$$ 
Then $nd^{2} \leq 8$, which is a contradiction. 
\epf

\subsection{Projective covers of simple left $D$-modules}
In this subsection we denote by $\widehat{D}$ the set of isomorphism classes 
of simple left $D$-modules and by $P(V)$ the projective cover of a simple $D$-module $V$.
The left regular $D$-module decomposes as 
$$_{D} D \simeq \bigoplus_{V\in \widehat{D}} P(V)^{\dim V}.
$$

\begin{lema}\label{lem:projcovers}
{\color{white} a}	
 \begin{enumerate} 
  \item[$(i)$] $V_{i,j}\ot \Bbbk_{\chi^{\ell}} 
  \simeq V_{i+\ell,j+ 2\ell}$ and 
  $\Bbbk_{\chi^{\ell}}\ot \Bbbk_{\chi^{k}} \simeq  \Bbbk_{\chi^{\ell+k}}$
  for all $(i,j)
  \in \Lambda$ and $k,\ell \in \Z_{4}$.
  \item[$(ii)$]
  $P(V_{i,j}) \simeq V_{i,j}$ 
  for all $(i,j)
  \in \Lambda$.
  \item[$(iii)$] $P(\Bbbk_{\chi^{\ell}}) \simeq  P(\Bbbk_{\eps})\ot \Bbbk_{\chi^{\ell}}$ and 
  $\dim P(\Bbbk_{\chi^{\ell}}) = 4$ for all $\ell\in \Z_{4}$.
 \end{enumerate}
\end{lema}

\pf  {\color{white} a}	

\noindent $(i)$ follows by a direct computation.

\noindent $(ii)$ Let $(i,j) \in \Lambda$ and $\mu=\chi^{\ell}$ be a character on $D$. Since
$\Hom (P(V_{i,j})\ot \Bbbk_{\mu}, V_{i,j}\ot\Bbbk_{\mu}) = 
\Hom (P(V_{i,j}), V_{i,j}\ot\Bbbk_{\mu}\ot \Bbbk_{\mu}^{*})
= \Hom (P(V_{i,j}), V_{i,j})\neq 0$, 
and $P(V_{i,j})\ot \Bbbk_{\mu}$ is projective, it follows that 
$P(V_{i,j})\ot \Bbbk_{\mu}$ contains  $P(V_{i,j}\ot\Bbbk_{\mu}) \simeq P(V_{i+\ell, j+2\ell})$.
This inclusion induces an isomorphism when tensoring with $\Bbbk_{\chi^{-\ell}}$. Hence,
$P(V_{i,j})\ot \Bbbk_{\mu}\simeq P(V_{i+\ell, j+2\ell})$.

Assume there is $(i,j) \in \Lambda$ such that  
$\dim P(V_{i,j}) > \dim V_{i,j}$. Since $D$ is unimodular, by \cite[page 487]{L} 
the socle of $P(V_{i,j})$ is
$V_{i,j}$.
Thus, $\dim P(V_{i,j})\geq 2 
\dim V_{i,j}$ and consequently
$ \dim P(V_{i-\ell,j-2\ell}) \geq 4$ for $\ell =0,1,2,3$.
Consider the set 
$$I = \{(m,n)\in \Lambda:\ (m,n) \neq 
(i+k,j+2k)\text{ for all }0\leq k\leq 3 \}.$$ 
It contains $8$ elements. We now have:
\begin{align*}
\dim D & = \sum_{j=0}^{3} \dim P(\Bbbk_{\chi^{j}}) + \sum_{(m,n) \in I}
2 \dim P(V_{m,n}) + 4\cdot 2 \dim P(V_{i,j})\\
& \geq   
\sum_{j=0}^{3} \dim P(\Bbbk_{\chi^{j}}) + 8\cdot 2\cdot 2 + 4\cdot 2\cdot 4 = \sum_{j=0}^{3} \dim P(\Bbbk_{\chi^{j}}) + 64,
\end{align*}
a contradiction. Hence, $\dim P(V_{i,j})= \dim V_{i,j} $ and consequently $P(V_{i,j})\simeq V_{i,j}$.

\noindent $(iii)$ The first assertion is well-known. 
Since $P(V_{i,j}) \simeq V_{i,j}$ 
  for all $(i,j)
  \in \Lambda$, we have that 
  $64 = 4\dim P(\Bbbk_{\eps}) + 48$, which implies that  $\dim P(\Bbbk_{\eps})=4$. 
\epf

We describe next the projective cover of the trivial module. Consider
 the $4$-dimensional $D$-module
 $P=\Bbbk\{p_{1},p_{2},p_{3},p_{4}\}$ whose structure is given by the following table:
\begin{equation}\label{eq:defproj}
\begin{tabular}{|c||c|c|c|c|}
\hline
$\cdot$ & $ p_{1} $ & $ p_{2} $ & $ p_{3} $ & $ p_{4} $ \\
\hline
\hline
$a$ & $ p_{1} $ & $-\xi p_{2} $ & $ \xi p_{3} $ & $ p_{4} $ \\
\hline
$ b $ & $ p_{3} $ & $ \xi p_{4} $ & $ 0 $ & $ 0 $  \\
\hline
$ c $ & $ -p_{3} $ & $ \xi p_{4} $ & $ 0 $ & $ 0 $  \\
\hline
$ d $ & $ p_{1} $ & $ \xi p_{2} $ & $ -\xi p_{3} $ & $ p_{4} $ \\
\hline
$ x $ & $ p_{2}+\sqrt{2}p_3 $ & $ -\sqrt{2}p_4 $ & $ p_{4} $ & $ 0 $ \\
\hline
$ g $ & $ p_{1} $ & $ -p_{2} $ & $ -p_{3} $ & $ p_{4} $ \\
\hline
\end{tabular}
 \end{equation}

Set $P_{i} = \Bbbk \{p_{i},\ldots, p_{4}\}$ for $i=1,2,3,4$.
It is easy to see that 
$P_{4}\simeq \Bbbk_{\eps}$, $P_{3}/P_{4} \simeq \Bbbk_{\chi}$,
$P_{2}/P_{3} \simeq \Bbbk_{\chi^{3}}$  and $P_{1}/P_{2} \simeq \Bbbk_{\eps}$ as $D$-modules.
In particular, $\Soc(P) \simeq   \Bbbk_{\eps} \simeq \Top(P)$.

\begin{lema}\label{lem:P-indec-series}
       $P$ is an indecomposable $D$-module with composition series given by 
       $P=P_{1} \supset  P_{2} \supset P_{3}\supset P_{4}\supset \{0\}$.
      \end{lema}

\pf Assume that $P=M\oplus N$ for two submodules $M$ and $N$ of $P$. Either $M$ or $N$
contains an element $w$ of the form 
$p_{1}+ \alpha p_{2} + \beta p_{3} + \gamma p_{4}$ with $\alpha, \beta, \gamma \in \Bbbk$. 
Without loss of generality, suppose that $w\in M$. Then 
$ p_4 =  (xb)\cdot w$ and $ p_3 + \xi\alpha p_4 =b\cdot w $ belong to $M$. This implies that 
$p_{3}, p_{4} \in M$ and, as a consequence, $p_{1} + \alpha p_{2} \in M$. The element 
$ x\cdot (p_1+\alpha p_2) = p_2 +\sqrt{2}p_{3} - \sqrt{2}\alpha p_4 $ belongs to  $M$ as well.
Hence, $p_{1}, p_{2}, p_{3}, p_{4} \in M$ and we are done.
The second assertion follows from the preceding discussion.
\epf

\begin{lema}
$P \simeq P(\Bbbk_{\eps})$ as $D$-modules.
\end{lema}

\pf
As $D$ is a Frobenius algebra, every projective module is injective; in particular, 
$P(\Bbbk_{\eps})$ is an injective envelope $E(\Bbbk_{\lambda})$ for some 
character $\lambda$. Moreover, since $D$ is unimodular, 
the socle and top of $P(\Bbbk_{\eps})$ coincide and we must have that 
$P(\Bbbk_{\eps}) \simeq E(\Bbbk_{\eps})$. On the other hand,
by Lemma \ref{lem:P-indec-series}, we know that  $P$ is an indecomposable module
with $\Soc(P) \simeq \Bbbk_{\eps}$. Thus, $P$ embeds in $E(\Bbbk_{\eps})$, which
implies that they are isomorphic, since they have the same dimension.
\epf

\begin{obs}\label{rmk:comp-fact-proj-cov}
Using that $P(\Bbbk_{\chi^{\ell}}) \simeq P\ot \Bbbk_{\chi^{\ell}}$ and  \eqref{eq:defproj}, one obtains a composition 
series of $P(\Bbbk_{\chi^{\ell}})$  by
tensoring the given composition series of $P$ with $\Bbbk_{\chi^{\ell}}$. Set $P_{j}(\Bbbk_{\chi^{\ell}}) = P_{j}\ot \Bbbk_{\chi^{\ell}}$
for all  $1\leq \ell \leq 3$ and $1\leq j\leq 4$. 
Then, one has that 
 $P_{3}(\Bbbk_{\chi^{\ell}})/P_{4}(\Bbbk_{\chi^{\ell}}) \simeq \Bbbk_{\chi^{\ell+1}}$,
$P_{2}(\Bbbk_{\chi^{\ell}})/P_{3}(\Bbbk_{\chi^{\ell}}) \simeq \Bbbk_{\chi^{\ell+3}}$  
and $P_{1}(\Bbbk_{\chi^{\ell}})/P_{2}(\Bbbk_{\chi^{\ell}}) \simeq \Bbbk_{\chi^{\ell}}$ as $D$-modules.
\end{obs}

As a direct consequence of the results above we have the following theorem.
\begin{teor}\label{thm:projectivecovers} The $D$-modules $P_{1}(\Bbbk_{\chi^{\ell}}) = P\ot \Bbbk_{\chi^{\ell}}$ and 
$V_{i,j}$, with $0\leq \ell\leq 3$ and $(i,j) \in \Lambda$, are the projective
covers of the simple $D$-modules. In particular,  
$$_{D} D \simeq \sum_{\ell=0}^{3} P_{1}(\Bbbk_{\chi^{\ell}}) \oplus 
\sum_{(i,j) \in \Lambda} V_{i,j}^{2}. 
$$
\end{teor}

\vspace{-0.8cm}\qed

\begin{obs}
For 
$0\leq \ell\leq 3$, let $\{p_{i,\ell}=p_{i}\ot 1\}_{1\leq i\leq 4}$ be
the linear basis of $P_{1}(\Bbbk_{\chi^{\ell}})$ constructed from the linear basis $\{p_{i}\}_{1\leq i\leq 4}$ of $P$. 
By \eqref{eq:defproj}, 
the $D$-module structure of $P_{1}(\Bbbk_{\chi^{\ell}})$ can be described explicitly:
\begin{equation}\label{eq:defprojall}
\begin{tabular}{l}
$ a\cdot (p_{i,\ell}) = (a\cdot p_i)\ot(a\cdot 1) + (c\cdot p_i)\ot(b\cdot 1) = \xi^\ell (a\cdot p_i)\ot 1$, \\
$ b\cdot (p_{i,\ell})  = (b\cdot p_i)\ot(a\cdot 1) + (d\cdot p_i)\ot(b\cdot 1) = \xi^{\ell} (b\cdot p_i)\ot 1 $, \\
$ c\cdot (p_{i,\ell}) = (a\cdot p_i)\ot(c\cdot 1) + (c\cdot p_i)\ot(d\cdot 1) = \xi^{-\ell} (c\cdot p_i)\ot 1$, \\
$ d\cdot (p_{i,\ell}) = (b\cdot p_i)\ot(c\cdot 1) + (d\cdot p_i)\ot(d\cdot 1) = \xi^{-\ell} (d\cdot p_i)\ot 1$, \\
$ x\cdot (p_{i,\ell})  = (g\cdot p_i)\ot(x\cdot 1) + (x\cdot p_i)\ot(1\cdot 1) =  (x\cdot p_i)\ot 1$, \\
$ g\cdot (p_{i,\ell})  = (g\cdot p_i)\ot(g\cdot 1) = (-1)^\ell (g\cdot p_i)\ot 1 $.
\end{tabular}
\end{equation}
\end{obs}

We end this subsection with the Clebsch-Gordan decomposition of the tensor product of two 2-dimensional simple $D$-modules.

\begin{prop}
 Let $V_{i,j}$ and $V_{k,l}$ be $2$-dimensional simple left $D$-modules. Then
  $$V_{i,j}\ot  V_{k,l}
  \simeq  \begin{cases}
  \begin{array}{ll}
    P(\Bbbk_{\chi^{i+k-1}}), & \text{ if } 2(i+k) +j+l \equiv 0\mod 4; \\
    V_{i+k,j+l} \oplus V_{i+k+3,j+l+2}, & \text{ otherwise.}
  \end{array}
 \end{cases}$$
  \end{prop}

\pf 
As $D$ is a 
quasitriangular Hopf algebra, by Lemma \ref{lem:projcovers} $(i)$ and Remark \ref{rmk:dualsVij},
we have that 
$\Hom_D (V_{i,j}\ot V_{k,l}, \Bbbk_{\chi^{t}})=
\Hom_D (V_{i,j}, \Bbbk_{\chi^{t}}\ot V_{k,l}^{*}) = \Hom_D (V_{i,j},  V_{-k+1+t,-l+2+2t})  $. 
By Schur's lemma, the latter is non-zero if and only if $2(i+k) + j+l = 0$ and $t=i+k-1$ in $\Z_{4}$. Also observe that $V_{i,j}\ot V_{k,l}$ is projective.

If $2(i+k) + j+l \equiv 0\mod 4$, then 
$\Hom (V_{i,j}\ot V_{k,l}, \Bbbk_{\chi^{i+k-1}}) \neq 0$ and, hence, $V_{i,j}\ot V_{k,l}$ must contain 
a submodule isomorphic to
$P(\Bbbk_{\chi^{i+k-1}})$. As both modules have the same dimension, the inclusion induces an
isomorphism.

Now, if $2(i+k) + j+l \not\equiv 0\mod 4$, then $V_{i,j}\ot V_{k,l}$ cannot contain a 1-dimensional module.
Thus, $V_{i,j}\ot V_{k,l}$ is isomorphic to a direct sum of
two 2-dimensional simple modules. Fix $\{v_{1},v_{2}\}$ and $\{w_{1},w_{2}\}$ linear bases of
 $V_{i,j}$ and $ V_{k,l} $, respectively, and set $u_{1}=v_{1}\ot w_{1}$,
 $u_{2}=v_{1}\ot w_{2}$, $u_{3}=v_{2}\ot w_{1}$ and $u_{4}=v_{2}\ot w_{2}$.
A direct computation shows that the matrices defining the actions $ \rho(a) $ and $ \rho(g) $ on $V_{i,j}\ot  V_{k,l}$, with respect to the
basis $\{u_{i}\}_{1\leq i \leq 4}$,
have the following form
 \[
\rho(a)= \left(
\begin{array}{cccc}
\xi^{i+k} & 0 &0&\xi^{2k} \\
0 & \xi^{i+k+3} &0&0\\
0 & 0& \xi^{i+k+3} &0\\
0 & 0&0& \xi^{i+k+2}\\
\end{array}
\right)
\text{ and } \ 
\rho(g)= \left(
\begin{array}{cccc}
\xi^{j+l} & 0 &0&0 \\
0 & -\xi^{j+l} &0&0\\
0 & 0& -\xi^{j+l} &0\\
0 & 0&0& \xi^{j+l}\\
\end{array}\right).
\]
Looking at the eigenspace decomposition with respect to the action of $a$ and $g$, it follows that necessarily 
$V_{i,j}\ot  V_{k,l}\simeq V_{i+k,j+l}\oplus V_{i+k+3,j+l+2}$. 
\epf

\subsection{Some indecomposable $D$-modules}
Let $A$ be a finite-dimensional $\Bbbk$-algebra and 
${V_{ 1} ,\ldots, V_{ n} }$ a complete list of non-isomorphic simple
left $A$-modules. The Ext-Quiver of $A$ is the quiver $\ExtQ(A)$
with vertices ${1,\ldots , n}$ and $\dim \Ext^{ 1}_{ A} (V_{ i} , V_{ j} )$ arrows from the vertex $i$ 
to the vertex $j$. 
Given a quiver $Q$ with vertices ${1,\ldots , n}$, its separation diagram is the
unoriented graph with vertices ${1,\ldots , n, 1' , \ldots , n' }$ and with an edge from $i$ to  
$j'$ for each arrow $i\to j$ in $Q$. The separation diagram of $A$
is the separation diagram of its Ext-Quiver.  It is well-known that 
a finite-dimensional algebra is of finite (tame) representation type if and only 
if its separation diagram
is a disjoint union of finite (affine) Dynkin diagrams.

In this section we compute the separation diagram of $D$ and show that $D$ is of tame
representation type. In order to do so, we use the isomorphism of abelian groups
between 
$\Ext^{ 1}_{ A} (V_{ i} , V_{ j} )$ and the
equivalence classes of extensions $0\to V_{j} \to M \to V_{i}\to 0$
of $V_{i}$ by $V_{ j}$.

\subsubsection{2-dimensional $($non-simple$)$ indecomposable modules}\label{subsec:2dimindec}
Let $A$ be the subalgebra of $D$
generated by $a,d$ and $g$. Then $A$ is an $8$-dimensional commutative
algebra given by $A=\Bbbk\langle a,g:\ a^{4}=1=g^{4}, ag=ga\rangle$.
In particular, all simple $A$-modules are $1$-dimensional.

\begin{defi}\label{def:2dimindec}
For  $0\leq \ell \leq 3$, let $M_{\ell}^{+}=\Bbbk\{m_{1},m_{2}\}$ be
the 2-dimensional $D$-module whose structure is given by setting $\Bbbk\, m_{1} \simeq \Bbbk_{\chi^{\ell}}$ and
\begin{align*}
a\cdot m_{2} & = \xi^{\ell+1}\ m_{2}, & b\cdot m_{2} & = 0= 
c\cdot m_{2},\\
g\cdot m_{2} & = (-1)^{\ell+1}\, m_{2}, &  
x\cdot m_{2} & = m_{1}.
\end{align*}
It is easy to see that $M_{\ell}^{+}$ is an indecomposable $D$-module that contains a submodule isomorphic to $\Bbbk_{\chi^{\ell}}$ and verifies that
$M_{\ell}^{+}/\Bbbk_{\chi^{\ell}} = \Bbbk_{\chi^{\ell+1}}$.
Analogously, for $0\leq \ell \leq 3$, let 
$M_{\ell}^{-}=\Bbbk\{m_{1},m_{2}\}$ be
the $2$-dimensional $D$-module given by $\Bbbk\, m_{1} \simeq \Bbbk_{\chi^{\ell}}$ and
\begin{align*}
a\cdot m_{2} &=  \xi^{\ell-1}\ m_{2}, & b\cdot m_{2} & 
= \frac{\sqrt{2}}{2} \xi^{\ell-1}\, m_{1}, &
c\cdot m_{2} &=
\frac{\sqrt{2}}{2}(- \xi)^{\ell+1} m_{1}, \\
g\cdot m_{2} & =(-1)^{\ell-1}\, m_{2},&
x\cdot m_{2} & =  m_{1}.
\end{align*} 
Then, $M_{\ell}^{-}$ is an indecomposable module that contains a submodule isomorphic to $\Bbbk_{\chi^{\ell}}$
and satisfies that $M_{\ell}^{-}/\Bbbk_{\chi^{\ell}} = \Bbbk_{\chi^{\ell-1}}$.
\end{defi}

Observe that, for all $0\leq \ell \leq 3$, the submodule $P_{3}(\Bbbk_{\chi^{\ell}})$ of $P_{1}(\Bbbk_{\chi^{\ell}})$ is isomorphic to $M^{+}_{\ell}$. 

\begin{lema} \label{lem:indecnonsimpledim2}
Fix $0\leq \ell \leq 3$.
\begin{enumerate}
 \item[$(i)$] Let $M$ be a 2-dimensional indecomposable $D$-module that contains a submodule isomorphic to $\Bbbk_{\chi^{\ell}}$. Then 
  $M\simeq M^{+}_{\ell}$
 or $M\simeq M^{-}_{\ell}$.
\item[$(ii)$] $$\dim \Ext^{ 1}_{D} ( \Bbbk_{\chi^{k}}, \Bbbk_{\chi^{\ell}} ) = \begin{cases}
\begin{array}{ll}
1, & \text{if } k=\ell \pm 1;\\
0, & \text{otherwise}.
\end{array}
\end{cases}
$$

 \end{enumerate}
\end{lema}

\pf 
$(i)$ Write $\lambda = \chi^{\ell}$. 
We must have that 
$M\simeq  \Bbbk_{\lambda} \oplus \Bbbk_{\mu}$ as $A$-modules, with $\mu$ some character on $A$. That is,  
$M$ has a linear basis $\{m_{1},m_{2}\}$ such that $\Bbbk m_{1} \simeq \Bbbk_{\lambda}$ (as $ D $-module) and 
$z\cdot m_{2} = \mu(z) \, m_{2}$, for any $z\in A$. 
Since $ b^2=0 $ and $ x^2=g^2-1 $, it follows that
\begin{align*}
b\cdot m_{2} &= \alpha\ m_{1}\qquad \text{ and } \qquad  x\cdot m_{2} = \beta\ m_{1},
\end{align*} 
for some $\alpha, \beta \in \Bbbk$. As $a^{2}b=c$, we have that $c\cdot m_{2} = \lambda(a)^{2} \alpha \, m_{1}$. 
Moreover, using the relation $bx -\xi xb = \sqrt{2}\xi(a - gd)$, it follows that $ \mu(g)=\mu(a)^2 $. Thus, $\mu$ 
defines a character $ \widetilde{\mu} $ on $D$ by taking $ \widetilde{\mu}(x) =\widetilde{\mu}(b) =\widetilde{\mu}(c) = 0 $ and $ \widetilde{\mu}|_A = \mu $.

Observe that $\alpha$ and $ \beta $ are not simultaneously zero. 
Indeed, if $ \alpha =\beta =0 $, then $M\simeq \Bbbk_{\lambda} \oplus \Bbbk_{\widetilde{\mu}}$ as $ D $-modules, 
a contradiction since $ M $ is indecomposable. 
Furthermore, the relation $ ax +\xi xa = \sqrt{2}\xi(b + gc) $ yields that 
\begin{align}\label{eq:dem_indecomp_mod_dim_2}
\beta(\lambda(a) + \xi \widetilde{\mu}(a))  = 2 \sqrt{2} \xi \alpha.
\end{align}
Hence, $\beta\neq 0$.

Moreover, from the relation $xg=-gx$, we get that $ \widetilde{\mu}(g) =-\lambda(g) $, since otherwise, we would get that $ \beta =0 $. Thus, $\widetilde{\mu}(a)^2 = (-1)^{\ell +1} $. This implies that 
$\widetilde{\mu}(a) = \pm \xi^{\ell+1}= \xi^{\ell\pm 1}$, and consequently $\widetilde{\mu} = \chi^{\ell\pm 1}$.

If $\widetilde{\mu} = \chi^{\ell+1}$, then  $\alpha =0 $, by \eqref{eq:dem_indecomp_mod_dim_2}. In this case, $M$ is an indecomposable module isomorphic to $M^{+}_{\ell}$. 
Denote this module by $M^{+}_{\ell}(\beta)$.

If $\widetilde{\mu} = \chi^{\ell-1}$, then $ \lambda(a) \beta = \sqrt{2}\xi \alpha$, by \eqref{eq:dem_indecomp_mod_dim_2}. 
Thus, $M$ is isomorphic to $M^{-}_{\ell}$ and $(i)$ follows. 
Denote this module by $M^{-}_{\ell}(\beta)$.

$(ii)$ By the preceding discussion, we have that  
$\dim \Ext^{1}_{D}(\Bbbk_{\chi^{k}},\Bbbk_{\chi^{\ell}})= 0 $ if $k\neq \ell \pm 1$. 
On the other hand, if
$M$ is a non-trivial extension of $\Bbbk_{\chi^{\ell}}$ by 
$\Bbbk_{\chi^{\ell \pm 1}}$, then 
$M= M^{\pm}_{\ell}(\beta)$ for 
some $\beta \in \Bbbk^{\times}$.  
Assume $ M^{\pm}_{\ell}(\beta)\simeq  M^{\pm}_{\ell}(\beta')$ as extensions,
with $\beta, \beta' \in \Bbbk^{\times}$. Let
$\{m_{1}, m_{2}\}$ and  $\{m_{1}', m_{2}'\}$ be the linear bases of 
$M^{\pm}_{\ell}(\beta)$ and $M^{\pm}_{\ell}(\beta')$, respectively,  as defined above; and write 
$\varphi: M^{\pm}_{\ell}(\beta) \to M^{\pm}_{\ell}(\beta')$ for the isomorphism. Then, 
we must have that $\varphi(m_{1}) = m_{1}'$ and
$\varphi(m_{2}) = \gamma\, m_{1}' + 
\eta\, m_{2}'$ for some $\eta\neq 0$. 
Moreover, since $\varphi(x\cdot m_{2}) = \beta \varphi(m_{1}) =\beta\, m_{1}'$
equals $x\cdot \varphi(m_{2}) = \eta\, \beta'\, m_{1}'$, it follows that
$\beta = \eta\beta' $. This implies that  
$\dim \Ext^{1}_{D}(\Bbbk_{\chi^{\ell\pm 1}},\Bbbk_{\chi^{\ell}})= 1$ and the lemma is proved.
\epf

\begin{lema}\label{lem:dimExtVijKl=0}
	\begin{enumerate}
		\item[$(i)$] $\dim \Ext^{1}_{D}(V_{i,j},V_{k,\ell}) =0 $ 
		for all $(i,j),(k,\ell) \in \Lambda$.
		\item[$(ii)$] $\dim \Ext^{1}_{D}(V_{i,j},\Bbbk_{\chi^{\ell}}) =0 =  \dim \Ext^{1}_{D}(\Bbbk_{\chi^{\ell}}, V_{i,j})$ 
		for all $(i,j) \in \Lambda$ and $\ell\in \Z_{4}$.
	\end{enumerate}
\end{lema}

\pf The proof follows easily from the fact that $V_{i,j}$ is projective 
for all $(i,j) \in \Lambda$.
\epf

\begin{cor}\label{cor:reptypeD}
	$D$ is of tame representation type. 
\end{cor}

\pf
Denote by $i$ the vertex corresponding to the character $\chi^{i}$ for all
$0\leq i \leq 3$.
Lemma \ref{lem:indecnonsimpledim2} implies that $\ExtQ(D)$ contains the quiver 
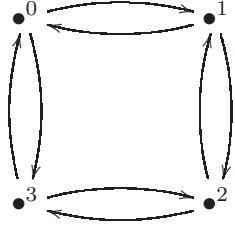
\begin{equation*}
\xymatrix{\bullet^0\ar@/^/[rr]\ar@/^/[dd]&
	&\bullet^1\ar@/^/[ll]\ar@/^/[dd]\\
	\\
	\bullet^3\ar@/^/[uu]\ar@/^/[rr]& &\bullet^2\ar@/^/[ll]\ar@/^/[uu] }
\end{equation*}
Thus, the separation diagram of $D$ contains the quiver $A^{(1)}_{3} \coprod A^{(1)}_{3}$
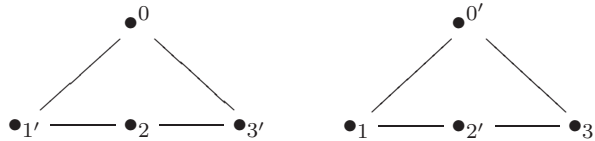
\begin{equation*}
\xymatrix{& \bullet^{0}\ar@{-}[ld]\ar@{-}[rd]&\\
	\bullet_{1'}\ar@{-}[r]& \bullet_{2}\ar@{-}[r]& \bullet_{3'}}\qquad
\xymatrix{& \bullet^{0'}\ar@{-}[ld]\ar@{-}[rd]&\\
	\bullet_{1}\ar@{-}[r]& \bullet_{2'}\ar@{-}[r]& \bullet_{3}}
\end{equation*}
Moreover, by Lemma \ref{lem:dimExtVijKl=0}, $\ExtQ(D)$ consists 
of the quiver above and twelve isolated points representing the simple modules 
$V_{i,j}$. Hence, $D$ is of tame representation type.
\epf

\subsubsection{3-dimensional indecomposable modules}\label{subsec:3dimindec}
In this subsection we describe the $3$-dimen\-sional indecomposable modules as we did
in the previous subsection for dimension $2$.

\begin{obs}\label{rmk:soctopone-dim}
Let $M$ be a non-simple 
indecomposable $D$-module.
As the simple $2$-dimensional $D$-modules  
$V_{i,j}$ are projective and injective
for all $(i,j) \in \Lambda$, they cannot be contained in any submodule of a quotient module of $M$.
In particular, $\Soc (M) $ and $\Top(M)$ consist of direct sums of 
$1$-dimensional modules.	

Furthermore, if $0\subseteq \Soc(M) \subseteq \Soc^{2}(M) \subseteq \cdots \subseteq \Soc^{n}(M) = M$
is the socle series of $M$,
then $\Soc(M/\Soc^{i}(M))$ does not contain any
 simple projective  and injective module. Indeed,
assume that $ \Soc(M/\Soc^{i}(M)) = \Soc^{i+1}(M)/\Soc^{i}(M)$ 
contains a simple projective and injective module $S$. Then, by the 
injectivity, $S$ is a direct summand of  $M/\Soc^{i}(M)$. Thus, it is a quotient
of $M/\Soc^{i}(M)$ and consequently also of $M$. Since $S$ is also projective, 
it must be a direct summand of $M$, a contradiction because $M$ is indecomposable.
\end{obs}

Next, we define some $3$-dimensional indecomposable $D$-modules.

\begin{defi}\label{def:3dimindec}
For $0\leq \ell \leq 3$, let
$N_{\ell}=\Bbbk\{n_{1},n_{2},n_{3}\}$ be
the 3-dimensional $D$-module whose structure is determined by setting $\Bbbk\, n_{1} \simeq \Bbbk_{\chi^{\ell}}$, $ \Bbbk\, n_{2}\simeq 
\Bbbk_{\chi^{\ell+2}}$ and
\begin{align*}
a\cdot n_{3} & = \xi^{\ell+1} \ n_{3}, & b\cdot n_{3} & = \frac{\sqrt{2}}{2}\xi^{\ell+1} n_{2},& 
c\cdot n_{3}& = -\frac{\sqrt{2}}{2}(-\xi)^{\ell+1} n_{2},\\ 
g\cdot n_{3} &= (-1)^{\ell+1}\, n_{3}, &  
x\cdot n_{3} &  = n_{1} + n_{2}.& & 
\end{align*}
It holds that $N_{\ell}$ is an indecomposable $D$-module with socle 
isomorphic to $\Bbbk_{\chi^{\ell}}\oplus \Bbbk_{\chi^{\ell+2}}$. Moreover, one has that 
$N_{\ell}/(\Bbbk_{\chi^{\ell}}\oplus \Bbbk_{\chi^{\ell+2}}) \simeq \Bbbk_{\chi^{\ell+1}}$,
$N_{\ell}/\Bbbk_{\chi^{\ell}} \simeq  M^{-}_{\ell + 2}$ and
$N_{\ell}/\Bbbk_{\chi^{\ell+2}} \simeq M^{+}_{\ell}$. 
\end{defi}

Recall from Remark \ref{rmk:comp-fact-proj-cov} that for 
$0\leq \ell \leq 3$, there is a 3-dimensional indecomposable $D$-module 
$P_{2}(\Bbbk_{\chi^{\ell}})$, which is the unique maximal submodule of
$P_{1}(\Bbbk_{\chi^{\ell}})$.

\begin{lema}\label{lem:3dimindec}
Let $N$ be a 3-dimensional indecomposable $D$-module. 
Then $N\simeq P_{2}(\Bbbk_{\chi^{\ell}})$ or $N\simeq N_{\ell}$ 
for some $0\leq \ell\leq 3$.
\end{lema}

\pf By Remark \ref{rmk:soctopone-dim}, $\Soc(N)$ contains 
only $1$-dimensional modules. 
Assume first that $\Soc(N) = \Bbbk_{\lambda}$ for some 
$D$-character $\lambda$. Then $N$ embeds in the injective hull $E(\Bbbk_{\lambda})$ of $\Bbbk_{\lambda}$, which 
is isomorphic to the projective cover $P_{1}(\Bbbk_{\lambda})$, because
$D$ is unimodular. Thus,
$N$ must be isomorphic to the unique maximal submodule $P_{2}(\Bbbk_{\lambda})$ of $P_{1}(\Bbbk_{\lambda})$.

Suppose now that $\Soc(N) = \Bbbk_{\lambda}\oplus \Bbbk_{\mu}$ for some 
$D$-characters $\lambda, \mu$. Then $N$ fits into an exact sequence 
\begin{equation}\label{eq:exseqindecN}
 0\to \Bbbk_{\lambda}\oplus \Bbbk_{\mu} \to N \to \Bbbk_{\tau} \to 0,
\end{equation}
for some 
$D$-character $\tau$. In particular, $N\simeq  \Bbbk_{\lambda}\oplus \Bbbk_{\mu} \oplus \Bbbk_{\tau}$
as $A$-modules.
Let $\{n_{1},n_{2},n_{3}\}$ be a linear basis of $N$ such that 
$\Bbbk n_{1} \simeq \Bbbk_{\lambda}$, $\Bbbk n_{2} \simeq \Bbbk_{\mu}$ and 
\begin{align*}
a\cdot n_{3} &= \tau(a) n_{3}, &d\cdot n_{3} &= \tau(d) n_{3},& 
g\cdot n_{3} &= \tau(g) n_{3}= \tau(a)^{2}n_{3},\\
b\cdot n_{3} &= \beta_{1} n_{1}+ \beta_{2} n_{2},&
c\cdot n_{3} &= \gamma_{1} n_{1}+ \gamma_{2} n_{2},&
x\cdot n_{3} &= \theta_{1} n_{1}+ \theta_{2} n_{2}.
\end{align*}
As $a^{2}b = c$, we have that $\lambda(a)^{2}\beta_{1} = \gamma_{1}$ and 
$\mu(a)^{2}\beta_{2} = \gamma_{2}$. Also,
using the relations $cg=-gc$, $bg=-gb$, $xg=-gx$ and $ax + \xi xa = \sqrt{2}\xi(b + gc)$,
we obtain the equalities 
\begin{align*}
\beta_{1}(\lambda(g) +\tau(g)) &= 0, &  \beta_{2}(\mu(g) +\tau(g)) & = 0,\\
\gamma_{1} (\lambda(g) +\tau(g)) & = 0, & \gamma_{2} (\mu(g) +\tau(g)) & = 0,\\
\theta_{1} (\lambda(g) +\tau(g) ) &=  0, & \theta_{2} (\mu(g) +\tau(g) ) &= 0,\\
 \theta_{1}(\lambda(a) +\xi \tau(a)) & = 2\sqrt{2}\xi \beta_{1}, & 
\theta_{2}(\mu(a) +\xi \tau(a))  & = 2\sqrt{2}\xi \beta_{2}.  
\end{align*}
If $\tau(g) \neq -\lambda(g)  $ and $\tau(g) \neq -\mu(g)  $, 
then $\beta_{i} = \gamma_{i} = \theta_{i} =0$
and consequently $N\simeq \Bbbk_{\lambda}\oplus \Bbbk_{\mu} \oplus \Bbbk_{\tau}$
as $D$-modules, a contradiction.
If $\tau(g) = -\lambda(g)  $ but $\tau(g) \neq -\mu(g)  $, then 
$\beta_{2} = \gamma_{2} = \theta_{2} =0$ and this implies that 
$N = L \oplus \Bbbk_{\mu}$ with $L = \Bbbk\{n_{1},n_{3}\}$. Analogously,
$N$ is decomposable if   $\tau(g) = -\mu(g)  $ but $\tau(g) \neq -\lambda(g)  $.
Hence, $-\tau(g) = \lambda(g) = \mu(g) $ and thus
$\lambda(a) = \pm \mu(a)= \pm \xi \tau(a)$ or $\lambda(a) = \pm \mu(a)= \mp \xi \tau(a)$. 
The same reasoning shows that 
$\theta_{1}\neq 0 \neq \theta_{2}$ since otherwise 
$N$ would be decomposable. So, we may assume that $\theta_{1} = \theta_{2} = 1$.

Moreover,  $\lambda = - \mu $ since otherwise 
$N$ is also decomposable. Indeed, if $\lambda = \mu$, we have that $ \beta_1 = \beta_2 $, from which follows that 
$N \simeq \Bbbk\{n_{1}\} \oplus \Bbbk\{v,n_{3}\}$ as $D$-modules with $v= n_{1}+ n_{2}$.

Set $\lambda = \chi^{\ell}$ for some 
$0 \leq\ell \leq 3 $. Then $\mu = \chi^{\ell+2}$ and $\tau = \chi^{\ell  \pm 3}= \chi^{\ell  \mp 1}$.
From the paragraph above it follows that $\beta_{1} = 0$ 
for $\tau(a) = \chi^{\ell+1}(a) = \xi^{\ell +1}$, and
$\beta_{2} = 0$ for $\tau(a) = \chi^{\ell-1}(a) = \xi^{\ell -1}$. If 
 $\tau=\chi^{\ell + 1}$, then $\beta_{2} = \frac{\sqrt{2}}{2}\xi^{\ell+1}$ and 
 $\gamma_{2} = -\frac{\sqrt{2}}{2}(-\xi)^{\ell+1}$. In such a case, $N\simeq N_{\ell}$.
 If $\tau=\chi^{\ell - 1}$, the same argument shows that 
$N\simeq N_{\ell +2}$ and the lemma is proved.
\epf

\begin{obs} Observe that $N_{\ell}^{*} \simeq P_{2}(\Bbbk_{\chi^{-\ell-1}})$ as $D$-modules,
since
$N_{\ell}^{*}$ is a 3-dimensional indecomposable module
with socle $\Bbbk_{\chi^{-\ell-1}}$.
\end{obs}

\section{The category $\ydk$}\label{sec:categoryyd}
Using 
the equivalence $_{D}\mathcal{M} \simeq \ydk$, we determine in this section the 
simple and some indecomposable objects of $\ydk$, and describe their braidings. Note that, 
by \cite[Proposition 2.2.1]{AG}, one has $\ydk \simeq \yda$ with 
$\mathcal{A} = \mathcal{A}''_{4}=\Kk^{*}$.

\subsection{Simple objects and projective covers in $\ydk$}
Our intention is to describe the simple $D$-modules and their projective covers  
as left Yetter-Drinfeld modules over $ \Kk $. To achieve our goal, we simply need to describe the coaction, since the action
is given by the restriction to $\Kk$ of the action of $D$.

\begin{prop}\label{prop:yddim 1}  For $0\leq j \leq 3$, set $ \Bbbk_{\chi^{j}} = \Bbbk{v}$. Then
$ \Bbbk_{\chi^{j}}\in \ydk$ with its structure given by 
$$a\cdot v = \xi^{j}v,\qquad b\cdot v = c\cdot v = 0,\qquad d\cdot v = \xi^{-j}v\quad\text{ and }\quad
\delta(v) =   
a^{2j}\otimes v.
$$
\end{prop}

\pf 
Since $\Bbbk_{\chi^{j}}$ is $1$-dimensional, we must have that $\delta(v) = h\ot v$ for some
$h\in G(\Kk) = \{1, a^{2}\}$. As $f\cdot v = \langle f,h\rangle v$ for all $f\in \Kk^{*}$ and 
$ \langle g,a^{2}\rangle = -1 $, it follows that 
$\delta(v)  =  a^{2j}\otimes v$.
\epf

\begin{prop} \label{prop:braided spaces of dim 1}
The braiding of $ \Bbbk_{\chi^{j}}$ is given by
$c(v\otimes v) = (-1)^{j}
v \otimes v
$.
\end{prop}

\pf Just apply formula \eqref{eq:braidingyd} to Proposition \ref{prop:yddim 1}.
\epf

Now we describe the Yetter-Drinfeld structure of the $2$-dimensional simple modules.

\begin{prop}\label{prop:yddim 2} For $(i,j) \in \Lambda$, set
$ V_{i,j}=\Bbbk\{e_{1},e_{2}\} $, and
write $\lambda_{1}=\xi^{i}$ and $\lambda_{2}=\xi^{j}$. Then $ V_{i,j} \in \ydk$ with its action given by
\begin{align*}
a\cdot e_{1} & = \lambda_1 e_{1}, &  
b\cdot e_{1} & = 0, & 
c\cdot e_{1} & = 0, & 
d\cdot e_{1} & = \lambda_1^{-1} e_{1}, \\
a\cdot e_{2} & =  -\xi\lambda_1 e_{2},  &
b\cdot e_{2} & =  \lambda_1^{2} e_{1},  & 
c\cdot e_{2} & =  e_{1}, & 
d\cdot e_{2} & = \xi \lambda_1^{-1} e_{2},
\end{align*}
and its coaction by
\begin{align*}
 \delta(e_{1}) &= 1\otimes e_{1} -2\lambda_{1} ac\otimes e_{2}, &  
\delta(e_{2}) &= a^2\otimes e_{2},& &\text{for }\lambda_{2} = 1, \\
 \delta(e_{1}) &= a^2\otimes e_{1} + 2\lambda_{1} ab\otimes e_{2}, &  
\delta(e_{2}) &= 1\otimes e_{2},& &\text{for }\lambda_{2} = -1, \\
 \delta(e_{1}) &= d\otimes e_{1} +(\lambda_{1}^3 -\xi\lambda_{1}) c\otimes e_{2}, &  
\delta(e_{2}) &= a\otimes e_{2} +\dfrac{1}{2}(\lambda_{1} +
\xi\lambda_{1}^3) b\otimes e_{1},& &\text{for }\lambda_{2} = \xi,\\
\delta(e_{1}) &= a\otimes e_{1} +(\lambda_{1}^3 +\xi\lambda_{1}) b\otimes e_{2}, &  
\delta(e_{2})& = d\otimes e_{2} +\dfrac{1}{2}(\lambda_{1} -\xi\lambda_{1}^3) c\otimes e_{1},
& &\text{for }\lambda_{2} = -\xi .\\
 \end{align*}
\end{prop}

\vspace{-0.7cm}

\pf 
Let $\{v_{i}\}_{1\leq i \leq 8}$ be a basis of $\Kk$ and 
$\{v^{i}\}_{1\leq i \leq 8 }$ its dual basis.
Recall that  
$\delta(v) = \sum_{i=1}^{8} v_{i}\ot v^{i}\cdot v$ for all $v\in V_{i,j}$.
Then, by the isomorphism from Lemma \ref{lem:isomAK}, we have that 
\begin{align}
\begin{split}\label{eq:coaction}
\delta(e_{1}) &= \sum_{i=0}^3 (g^{i})^{*} \otimes g^i\cdot e_{1}  + 
\sum_{i=0}^3 (xg^{i})^{*} \otimes xg^i \cdot e_{1} \\
&= \sum_{i=0}^3 (g^{i})^{*} 
\otimes \lambda_{2}^i e_{1} + \sum_{i=0}^3 (xg^{i})^{*} \otimes (-\lambda_{2})^i x_{21} e_{2}, \\
\delta(e_{2}) &= \sum_{i=0}^3 (g^{i})^{*} \otimes  g^i \cdot e_{2} + 
\sum_{i=0}^3 (xg^{i})^{*} \otimes xg^i\cdot e_{2}\\
& = \sum_{i=0}^3(g^{i})^{*} 
\otimes (-\lambda_{2})^i e_{2} + \sum_{i=0}^3 (xg^{i})^{*} \otimes \lambda_{2}^i x_{12} e_{1},
\end{split}
\end{align}
where $ (g^{i})^{*} = \dfrac{1}{4}(1 + \xi^i a + (-\xi)^i d + (-1)^i a^2) $, 
$ (xg^{i})^{*} = \dfrac{1}{4\sqrt{2}\xi}((-\xi)^i b + \xi^i c + ab + (-1)^i ac)$ for all 
$0\leq i \leq 3$, and 
$ x_{21} = \sqrt{2}\xi(\lambda_{1}^3 -\lambda_{1}\lambda_{2}) $,
$ x_{12} = \dfrac{\sqrt{2}}{2}\xi(\lambda_{1} +\lambda_{1}^{3}\lambda_{2})$.

Computing the formulae \eqref{eq:coaction} for each $ 0\leq i,j \leq 3 $, 
we obtain the explicit coactions presented in the statement of the proposition.
\epf

Next, we describe the braiding of the simple modules 
$V_{i,j}$ in $\ydk$. We use a 
matrix-like notation to state it in a compact form. Its proof follows by a direct computation 
using formula \eqref{eq:braidingyd} and Proposition \ref{prop:yddim 2}.

\begin{prop} \label{prop:braidedspacesdim2}
Set $\lambda_{1}=\xi^{i}$ and $\lambda_{2}=\xi^{j}$. The braiding of $V_{i,j} \in \ydk$ 
is given by the following formulae:
\begin{enumerate}
\item[$(i)$] For  $j=0$, $i \in \{1,3\} $: 
\[c(\left\{
\begin{array}{c}
e_1 \\
e_2 \\
\end{array}
\right\}\otimes \left\{
\begin{array}{cc}
e_1 & e_2 \\
\end{array}
\right\}) = \left\{
\begin{array}{cc}
e_1\otimes e_1 & e_2\otimes e_1 +2e_1\otimes e_2 \\
-e_1\otimes e_2 & e_2\otimes e_2 \\
\end{array}
\right\}.
\]
\item[$(ii)$] For $j=2$, $i \in \{0,2\} $: 
\[c(\left\{
\begin{array}{c}
e_1 \\
e_2 \\
\end{array}
\right\}\otimes \left\{
\begin{array}{cc}
e_1 & e_2 \\
\end{array}
\right\}) = \left\{
\begin{array}{cc}
e_1\otimes e_1 & -e_2\otimes e_1 +2e_1\otimes e_2 \\
e_1\otimes e_2 & e_2\otimes e_2 \\
\end{array}
\right\}.
\]
\item[$(iii)$] For $j=1$ and $i$ arbitrary: 
\[c(\left\{
\begin{array}{c}
e_1 \\
e_2 \\
\end{array}
\right\}\otimes \left\{
\begin{array}{cc}
e_1 & e_2 \\
\end{array}
\right\}) = \left\{
\begin{array}{cc}
\lambda_1^{3} e_1\otimes e_1 & \xi\lambda_1^{3}e_2\otimes e_1 +(\lambda_1^{3}-\xi\lambda_1) e_1\otimes e_2 \\
\lambda_1 e_1\otimes e_2 & -\xi\lambda_1 e_2\otimes e_2 +\dfrac{1}{2}(\lambda_1^{3}+\xi\lambda_1) e_1\otimes e_1 \\
\end{array}
\right\}.
\]
\item[$(iv)$] For $j=3$ and $i$  arbitrary:  
\[c(\left\{
\begin{array}{c}
e_1 \\
e_2 \\
\end{array}
\right\}\otimes \left\{
\begin{array}{cc}
e_1 & e_2 \\
\end{array}
\right\}) = \left\{
\begin{array}{cc}
\lambda_1 e_1\otimes e_1 & -\xi\lambda_1 e_2\otimes e_1 +(\lambda_1+\xi\lambda_1^{3}) e_1\otimes e_2 \\
\lambda_1^3 e_1\otimes e_2 & \xi\lambda_1^3 e_2\otimes e_2 +\dfrac{1}{2}(\lambda_1 -\xi\lambda_1^3) e_1\otimes e_1 \\
\end{array}
\right\}.
\]
\end{enumerate}\qed
\end{prop}

\begin{obs} \label{rmk:braidingsnotdiagonaltype}
All the braidings given by Proposition \ref{prop:braidedspacesdim2} are not of diagonal type. 
See the Appendix of the first \texttt{arXiv} version of this paper.
\end{obs}

We end this section with the description of the projective covers of the $1$-dimensional 
modules as objects in $\ydk$ and its braidings. 

Recall that the module 
$P_{1}(\Bbbk_{\chi^{j}}) $ is isomorphic to $P(\Bbbk_{\chi^{j}})$, for $0\leq j \leq 3$.

\begin{prop}\label{prop:ydprojcov}  
$ P_{1}(\Bbbk_{\chi^{j}}) \in \ydk$ with its action given by \eqref{eq:defproj}, \eqref{eq:defprojall};  and 
its coaction by 
\begin{align*}
\delta(p_{1,j}) & = (a^2)^j\ot p_{1,j} -\dfrac{\xi\sqrt{2}}{2}(a^2)^jac\ot (p_{2,j} + \sqrt{2}p_{3,j}),\\
\delta(p_{2,j}) & = (a^2)^{j+1}\ot p_{2,j} +\xi(a^2)^jab\ot p_{4,j},\\
\delta(p_{3,j}) & = (a^2)^{j+1}\ot p_{3,j} -\dfrac{\xi\sqrt{2}}{2}(a^2)^jab\ot p_{4,j}, \\
\delta(p_{4,j}) &= (a^2)^{j}\ot p_{4,j}.
\end{align*}
\end{prop}

\pf By the same reason presented in the beginning of the proof of the Proposition \ref{prop:yddim 2}, we obtain that
\begin{align*}
\delta(p_{1,j}) &= \sum_{i=0}^3 (g^{i})^{*} \otimes g^i \cdot p_{1,j}  +
\sum_{i=0}^3 (xg^{i})^{*} \otimes xg^i \cdot p_{1,j}   \\
&= \sum_{i=0}^3 (g^{i})^{*}\otimes ((-1)^j)^i p_{1,j} + 
\sum_{i=0}^3 (xg^{i})^{*} \otimes ((-1)^{j+1})^i  (p_{2,j}+\sqrt{2}p_{3,j}) \\
&= (a^2)^j\ot p_{1,j} -\dfrac{\xi\sqrt{2}}{2}(a^2)^jac\ot (p_{2,j} + \sqrt{2}p_{3,j}); \\
\delta(p_{2,j}) &= \sum_{i=0}^3 (g^{i})^{*} \otimes g^i \cdot	p_{2,j}  +
\sum_{i=0}^3 (xg^{i})^{*} \otimes   xg^i \cdot p_{2,j}  \\
&= \sum_{i=0}^3 (g^{i})^{*}\otimes ((-1)^{j+1})^i p_{2,j} + 
\sum_{i=0}^3 (xg^{i})^{*} \otimes (-\sqrt{2}((-1)^j)^i)p_{4,j} \\
&= (a^2)^{j+1}\ot p_{2,j} +\xi(a^2)^jab\ot p_{4,j}; \\
\delta(p_{3,j}) &= \sum_{i=0}^3 (g^{i})^{*} \otimes g^i  \cdot p_{3,j} +
\sum_{i=0}^3 (xg^{i})^{*} \otimes xg^i \cdot  p_{3,j}  \\
&= \sum_{i=0}^3 (g^{i})^{*}\otimes ((-1)^{j+1})^i p_{3,j} + 
\sum_{i=0}^3 (xg^{i})^{*} \otimes ((-1)^j)^ip_{4,j} \\
&= (a^2)^{j+1}\ot p_{3,j} -\dfrac{\xi\sqrt{2}}{2}(a^2)^jab\ot p_{4,j}; \\
\delta(p_{4,j}) &= \sum_{i=0}^3 (g^{i})^{*} \otimes g^i \cdot p_{4,j}  
+\sum_{i=0}^3 (xg^{i})^{*} \otimes xg^i \cdot p_{4,j}   = 
\sum_{i=0}^3 (g^{i})^{*}\otimes ((-1)^{j})^i p_{4,j}\\
&= (a^2)^{j}\ot p_{4,j}.
\end{align*}
\epf

The following result holds by a straightforward computation using \eqref{eq:braidingyd}
and the coactions given in Proposition \ref{prop:ydprojcov}.

\begin{prop} \label{prop:braidedspacesprojdim4} Fix $0\leq j\leq 3$.
The braiding of $P_{1}(\Bbbk_{\chi^{j}})\in \ydk$ 
is given by the formulae:
\begin{align*}
c(p_{1,j}\ot\left\{
\begin{array}{c}
p_{1,j} \\
p_{2,j} \\
p_{3,j} \\
p_{4,j} \\
\end{array}
\right\}) &= \left\{
\begin{array}{c}
(-1)^{j} p_{1,j} \\
p_{2,j} \\
p_{3,j} \\
(-1)^{j}p_{4,j} \\
\end{array}
\right\}\ot p_{1,j} +\dfrac{\sqrt{2}}{2} \left\{
\begin{array}{c}
-p_{3,j} \\
(-1)^{j}p_{4,j} \\
0 \\
0 \\
\end{array}
\right\}\ot(p_{2,j}+\sqrt{2}p_{3,j}),\\
c(p_{2,j}\ot\left\{
\begin{array}{c}
p_{1,j} \\
p_{2,j} \\
p_{3,j} \\
p_{4,j} \\
\end{array}
\right\}) 
&= \left\{
\begin{array}{c}
p_{1,j} \\
(-1)^{j+1}p_{2,j} \\
(-1)^{j+1}p_{3,j} \\
p_{4,j} \\
\end{array}
\right\}\ot p_{2,j} +\left\{
\begin{array}{c}
(-1)^{j+1}p_{3,j} \\
-p_{4,j} \\
0 \\
0 \\
\end{array}
\right\}\ot p_{4,j},\\
c(p_{3,j}\ot\left\{
\begin{array}{c}
p_{1,j} \\
p_{2,j} \\
p_{3,j} \\
p_{4,j} \\
\end{array}
\right\}) &= \left\{
\begin{array}{c}
p_{1,j} \\
(-1)^{j+1}p_{2,j} \\
(-1)^{j+1}p_{3,j} \\
p_{4,j} \\
\end{array}
\right\}\ot p_{3,j} +\dfrac{\sqrt{2}}{2}\left\{
\begin{array}{c}
(-1)^{j}p_{3,j} \\
p_{4,j} \\
0 \\
0 \\
\end{array}
\right\}\ot p_{4,j},\\
c(p_{4,j}\ot\left\{
\begin{array}{c}
p_{1,j} \\
p_{2,j} \\
p_{3,j} \\
p_{4,j} \\
\end{array}
\right\}) & =\left\{
\begin{array}{c}
(-1)^{j}p_{1,j} \\
p_{2,j} \\
p_{3,j} \\
(-1)^j p_{4,j}\\
\end{array}
\right\}\ot p_{4,j}.
\end{align*}\qed 
\end{prop}

\section{Nichols algebras in $\ydk$}\label{sec:NicholsalgK}
In this section we determine all finite-dimensional Nichols algebras of simple 
modules over $\Kk$. They consist of exterior 
algebras of dimension 2 and 8-dimensional algebras with triangular braiding. 
Indeed,
since all objects 
in $\ydk$ can be described as objects in the category of Yetter-Drinfeld modules over
the pointed Hopf algebra
$\Kk^{*}=\mathcal{A}''_{4}$, by \cite{ufer} it follows that the associated braiding is triangular. To the best of our knowledge, 
these 8-dimensional Nichols algebras constitute new examples. They 
are isomorphic to quantum linear spaces as algebras, but not as coalgebras since the
braiding differs; in our case, the braiding is not diagonal. 
It remains an open question if they are twist equivalent and in such a case, in which category.

We begin by studying the Nichols algebras of the $1$-dimensional simple modules
and their projective covers.

\begin{lema}\label{lema:nichols-1-dim} Let $0\leq j \leq 3$. 
The Nichols algebras $ \mathfrak{B}(\Bbbk_{\chi^{j}}) $ associated with $\Bbbk_{\chi^{j}}=\Bbbk x$ are:
\[ \mathfrak{B}(\Bbbk_{\chi^{j}}) = \left\{
\begin{array}{ll}
\Bbbk[x], & \text{for }j = 0, 2, \\
\Bbbk[x]/(x^{2})= \bigwedge \Bbbk_{\chi^{j}}, &\text{for } j = 1, 3.\\
\end{array}
\right.
\]
\end{lema}

\pf Immediate, since the braiding $ c = (-1)^j\tau$, where 
$ \tau $ represents the usual flip.
\epf

\begin{cor}\label{cor:prodnicholsdim1}
 Let $W= \Bbbk_{\chi^{j_{1}}}\oplus \cdots \oplus \Bbbk_{\chi^{j_{t}}}$ be 
 a direct sum of $1$-dimensional modules with $j_{s} \in \{1,3\}$ for all $1\leq s\leq t$. 
 Then $\toba(W) = \bigwedge W \simeq \toba(\Bbbk_{\chi^{j_{1}}}) \underline{\ot} \cdots \underline{\ot} \toba(\Bbbk_{\chi^{j_{t}}})  $.
\end{cor}

\pf If $j_{s} \in \{1,3\}$ for all $1\leq s\leq t$,
then the braiding on $W\ot W$ is $-\tau$ and therefore $\toba(W) = \bigwedge W  $. 
The last assertion follows from \cite[Theorem 2.2]{Gr}, because $c_{W\ot W}^{2}=\id_{W\ot W}$. Indeed,
if $v \in  \Bbbk_{\chi^{j_{r}}}$ and $w \in  \Bbbk_{\chi^{j_{s}}}$, then $c(v\ot w) = (a^{2})^{j_{r}}w\ot v =
(-1)^{j_{r}j_{s}}w\ot v$. 
\epf

\begin{lema}\label{lema:nichols-proj-1-dim}  
Let $0\leq j \leq 3$. Then 
$\mathfrak{B}(P(\Bbbk_{\chi^{j}})) $ is infinite-dimensional.
\end{lema}

\pf
In all cases, the braiding on $P(\Bbbk_{\chi^{j}})\ot P(\Bbbk_{\chi^{j}})$ contains an eigenvector of eigenvalue
$1$. The claim then follows by Remark \ref{def:nicholsbraiding}. 
Indeed, by Proposition \ref{prop:braidedspacesprojdim4} we have
that $c(p_{4,j}\ot p_{4,j}) = p_{4,j}\ot p_{4,j}$ for $j=0,2$, and 
$c(p_{3,j}\ot p_{3,j}) = p_{3,j}\ot p_{3,j}$ for $j=1,3$. 
\epf

Before we describe the Nichols algebras associated with $2$-dimensional simple modules, we analyze 
the Nichols algebras of non-simple indecomposable modules. It turns out that they are all infinite-dimensional.

\begin{obs}\label{rmk:soctopindecnichols}
 Let $V \in \ydk$ be a finite-dimensional module such that $\dim \toba(V)<\infty$.
 Since taking the Nichols algebra defines a functor between the category of braided vector spaces and the 
 category of braided Hopf algebras, see \cite{AG}, it follows that $\dim \toba(W) <\infty$ 
 for all $W \in \Soc(V)$ or $W \in \Top(V)$.
 Furthermore, let $0\subseteq \Soc(V) \subseteq \Soc^{2}(V) \subseteq \cdots \subseteq \Soc^{n}(V) = V$
 be a socle series of $V$. Then $\dim \toba(V/\Soc^{i}(V)) $, $\dim \toba(\Soc^{i}(V))$ 
 and $\dim \toba(\Soc(V/\Soc^{i}(V))) $ are finite 
 for all $1\leq i \leq n$.
\end{obs}

\begin{teor}\label{thm:indecnicholsinfdim}
 Let $M \in \ydk$ be a finite-dimensional non-simple indecomposable module. 
 Then $\toba(M)$ is infinite-dimensional.
\end{teor}

\pf We prove the claim by induction on $\dim M$.
Assume first that $\dim M = 2$. By Lemma \ref{lem:indecnonsimpledim2} $(i)$, we have that 
$M\simeq M^{+}_{\ell}$ or $M\simeq M^{-}_{\ell}$ for some $0\leq \ell \leq 3$. 
Since $\Soc(M^{+}_{\ell}) = \Bbbk_{\chi^{\ell}}$,
$\Top(M^{+}_{\ell}) = \Bbbk_{\chi^{\ell+1}}$ and $\Soc(M^{-}_{\ell}) = \Bbbk_{\chi^{\ell}}$,
$\Top(M^{-}_{\ell}) = \Bbbk_{\chi^{\ell-1}}$, by Lemma \ref{lema:nichols-1-dim} and Remark \ref{rmk:soctopindecnichols},
it follows that $\toba(M)$ is infinite-dimensional.


Assume now that $\dim M = n \geq  3$ and 
suppose that $\dim \toba(N)$ is infinite for all
indecomposable module of dimension less than $n$.  
By Remark \ref{rmk:soctopone-dim},
$\Soc(M)$ consists of $1$-dimensional modules. Let $\bar{N}$ be a simple module 
contained in $\Soc(M/\Soc(M))$ and denote by $N$ the corresponding submodule 
of $M$. Also $\dim \bar{N} = 1$ by Remark \ref{rmk:soctopone-dim}. 
If $\Soc(M)=\Bbbk_{\lambda}$, then 
$N$ is an indecomposable module of dimension $2$. The previous paragraph implies that
$\dim \toba(N)$ is infinite and consequently $\dim \toba(M)$ is infinite. 
Assume that $\Soc(M)$ contains more than one simple module and let  
$\Bbbk_{\lambda} \subset \Soc(M)$. If
$N/\Bbbk_{\lambda}$ is semisimple, then $N$ contains an indecomposable module
of dimension $2$ and whence $\dim \toba(N/\Bbbk_{\lambda})$ is infinite. 
This implies again that $\dim \toba(N)$ and $\dim \toba(M)$ are both infinite. 
If 
$N/\Bbbk_{\lambda}$ is not semisimple, then it contains an indecomposable module
of dimension less than $n$. By induction, $\dim \toba(N/\Bbbk_{\lambda})$ is infinite and 
the theorem follows.
\epf

\begin{obs}\label{rmk:nicholsss}
Let $V \in \ydk$ be such that $\dim \toba(V) $ is finite. Then by Theorem \ref{thm:indecnicholsinfdim}, 
$V$ is necessarily semisimple. In these notes we analyse only Nichols algebras 
over simple modules, since the case of semisimple modules demands much more
work to be carried out. A first approach could be done by studying the Yetter-Drinfeld submodules 
$\ad^{n}(V)(W)$ of a given Nichols algebra $\toba(V\oplus W)$ with $V$ and $W$ simple modules,
see \cite{HS} for details. A direct computation shows that $\toba(V\oplus W)$ is infinite-dimensional for 
$V = \Bbbk_{\chi}$, $W= V_{3,1}$, $V_{3,3}$, and 
$V = \Bbbk_{\chi^{3}}$, $W= V_{2,1}$, $V_{2,3}$. In fact, 
$\ad(\Bbbk_{\chi}) (V_{3,1}) \simeq V_{0,3}$, $\ad(\Bbbk_{\chi}) (V_{3,3}) \simeq V_{0,1}$,
$\ad(\Bbbk_{\chi^{3}}) (V_{2,1}) \simeq V_{1,3}$ and $\ad(\Bbbk_{\chi^{3}}) (V_{2,3}) \simeq V_{1,1}$.
\end{obs}

Now, we analyze the Nichols algebras associated with $2$-dimensional 
simple modules.  

\begin{lema}\label{lema:2dim-infinitenichols}
Let $\Lambda' = \Lambda \smallsetminus \{(2, 1), (3, 1),  (2, 3), (3, 3) \}$. 
Then $\toba(V_{i,j})$ is infinite-dimensional for all 
$(i,j)\in \Lambda'$.
\end{lema}

\pf In all cases, the braiding of $V_{i,j}$ contains an eigenvector $w\ot w$ of
eigenvalue 1, hence the lemma follows by Remark \ref{def:nicholsbraiding}. 
Indeed, for $ (i,j) = (1, 1) $ or $(1,3)$, the element
$w=e_1 +\sqrt{2}\xi e_2 $ do the job. For the other cases, take $w=e_{1}$.
\epf

Next, we describe the Nichols algebras associated with the pairs in $\Lambda \backslash \Lambda' $ by generators and relations.
It turns out that all of them are isomorphic to algebras associated with quantum linear planes.

Recall that every graded Hopf algebra in $\ydk$ satisfies the Poincar\'e duality \cite[Proposition 3.2.2]{AG}, 
that is,  for $ R= \bigoplus_{i=0}^{N}R^{i}$ with 
$R^{N}\neq\{0\}$, it holds that 
$\dim R^{i} = \dim R^{N-i}$. 

\begin{prop}\label{prop:BV21}
$\mathfrak{B}(V_{2, 1})$ is the algebra generated by the elements $x,y$ satisfying the 
following relations 
\begin{equation}\label{eq:relBV21}
 x^2 = 0,\qquad xy +\xi yx = 0, \qquad y^4 = 0.
\end{equation}
In particular $ \dim\toba(V_{2, 1}) = 8 $.
\end{prop}

\pf Write $x=e_1$, $y=e_2$ for the linear generators of $ V_{2, 1} $. 
Then, by Proposition \ref{prop:braidedspacesdim2}, we have that
\[c(\left\{
\begin{array}{c}
x \\
y \\
\end{array}
\right\}\otimes \left\{
\begin{array}{cc}
x & y \\
\end{array}
\right\}) = \left\{
\begin{array}{cc}
- x\otimes x & -\xi y\otimes x +(\xi -1) x\otimes y \\
- x\otimes y & \xi y\otimes y -\dfrac{1}{2}(1+\xi) x\otimes x \\
\end{array} 
\right\}.
\]
Hence, the relations \eqref{eq:relBV21} must 
hold in $ \mathfrak{B}(V_{2,1}) $. Indeed,
the first two ones are easily checked since they are primitive elements of degree 2. 
Let us focus in the last one; we show that it is also primitive, modulo relations in degree $ 2 $. Since
\begin{align*}
\com(y^2) &=(y\ot 1+ 1\ot y)(y\ot 1 + 1\ot y)  = y^2\ot 1 + (1 + \xi) y\ot y -\dfrac{1}{2}(1 + \xi) x\ot x + 1\ot y^2,
\end{align*}
we get that 
\begin{align*}
\com(y^3) &= (y\ot 1 + 1\ot y)(y^2\ot 1 + (1 + \xi) y\ot y -
\dfrac{1}{2}(1 + \xi) x\ot x + 1\ot y^2) \\
 &= y^3\ot 1 +\dfrac{1}{2}(1+\xi) xy\ot x +\xi y^2\ot y -
 \dfrac{1}{2}(1+\xi) x\ot xy +\xi y\ot y^2 + 1\ot y^3 
\end{align*}
because $ c(y\ot y^2) = (x^2-y^2)\ot y + \dfrac{1}{2}(1-\xi)(yx - xy)\ot x $. Thus, as $ c(y\ot xy) = -\xi xy\ot y + \dfrac{1}{2}(1+\xi)x^2\ot x $ and $ c(y\ot y^3) = \xi(-y^3 + x^2y+yx^2-xyx)\ot y +
\dfrac{1}{2}(1+\xi)(y^2x -x^3 - yxy + xy^2)\ot x $, it follows that
\begin{align*} 
\com(y^4) &= (y\ot 1+ 1\ot y)\com(y^3) = y^4\ot 1 + 1\ot y^4.
\end{align*}

Hence, we have a graded braided Hopf algebra epimorphism 
$ \pi : T(V_{2,1}) / I \twoheadrightarrow \mathfrak{B}(V_{2,1}) $, where 
$ I $ is the two-sided ideal generated by the relations \eqref{eq:relBV21}.
Set $R=T(V_{2,1})/I$, then clearly $ R= \bigoplus_{i=0}^{4}R^{i}$ with $R^{4}\neq 0$, $R^{0}\simeq\Bbbk$ and  
$R^{1}\simeq V_{2,1}$.
By the Poincar\'e duality, we have that $\pi$ is injective in degree $ 0 $, $ 1 $, $ 3 $ and $ 4 $.
In order to prove that $ \pi $ is an isomorphism, it remains to show that $ \pi $ is injective in degree $ 2 $. 
This is equivalent to check that the relations in degree $ 2 $ in the Nichols algebra are just 
$ x^2 = 0 $ and $ xy +\xi yx = 0 $, which follows by a direct computation using the braiding.
\epf

The proof of the next three propositions follows the same lines of Proposition \ref{prop:BV21}. 
Thus, for their proof we only show that the defining relations hold in the Nichols algebra.
We also write $x=e_1$, $y=e_2$ for the linear generators of each $ V_{i, j} $ in its respective proof.

\begin{prop}\label{prop:BV23}
$\toba(V_{2, 3})$ is the algebra generated by the elements $x,y$ satisfying the following relations:
$$
x^2 = 0,\qquad xy -\xi yx = 0,\qquad y^4 = 0.
$$
In particular,
$\dim \toba(V_{2, 3}) = 8 $.
\end{prop}

\pf 
Using the braiding given by Proposition \ref{prop:braidedspacesdim2}, we see that
\begin{align*}
\com(x^2) &= x^2\ot 1 + 1\ot x^2, \\
\com(xy) &= xy\ot 1 -\xi x\ot y +\xi y\ot x + 1\ot xy, \\
\com(yx) &= yx\ot 1 + y\ot x - x\ot y + 1\ot yx, \\
\com(y^2) &= y^2\ot 1 + (1 -\xi) y\ot y +\dfrac{1}{2}(\xi - 1) x\ot x + 1\ot y^2.
\end{align*}
Thus, the relations $ x^2 = 0 $ and $ xy -\xi yx = 0 $ must hold in $ \mathfrak{B}(V_{2,3}) $, since both elements are
primitive of degree 2.
Let us check that the relation $y^{4}=0$ also holds.
Since
$ c(y\ot y^2) = (x^2-y^2)\ot y + \dfrac{1}{2}(1+\xi)(yx - xy)\ot x $, we have that
\begin{align*}
\com(y^3) &= y^3\ot 1 +\dfrac{1}{2}(1-\xi) xy\ot x -\xi y^2\ot y +\dfrac{1}{2}(\xi -1) x\ot xy -\xi y\ot y^2 + 1\ot y^3.
\end{align*}
Then, 
\begin{align*}
\com(y^4) &= \com(y^{3})(y\ot 1+ 1\ot y) =  \ y^4\ot 1 + 1\ot y^4,
\end{align*}
because $ c(y\ot y^3) = -\xi(-y^3 + x^2y+yx^2-xyx)\ot y + \dfrac{1}{2}(1-\xi)(y^2x -x^3 - yxy + xy^2)\ot x $ 
and $ c(y\ot xy) = \xi xy\ot y + \dfrac{1}{2}(1-\xi)x^2\ot x $. 
Hence, the relation $ y^4 = 0 $ must hold in $ \mathfrak{B}(V_{2,3}) $. 
\epf

\begin{prop}\label{prop:BV31} 
$\toba(V_{3,1})$ is the algebra generated by the elements $x, y$ 
satisfying the following relations: 
$$x^2 -2y^2 = 0, \qquad xy + yx = 0,\qquad y^4 = 0.$$
In particular, $\dim \toba(V_{3,1})  = 8 $.
\end{prop}
\pf  
By Proposition \ref{prop:braidedspacesdim2}, we get that 
\begin{align*}
\com(x^2) &= x^2\ot 1 +(1+\xi) x\ot x + 1\ot x^2,&
\com(xy) &= xy\ot 1 +\xi x\ot y - y\ot x + 1\ot xy, \\
\com(yx) &= yx\ot 1 + y\ot x -\xi x\ot y + 1\ot yx,&
\com(y^2) &= y^2\ot 1 +\dfrac{1}{2}(1 + \xi) x\ot x + 1\ot y^2.
\end{align*}
From this formulae it follows that the relations $ x^2 -2y^2 = 0 $ and $ xy + yx = 0 $ hold in $\toba(V_{3,1}) $.
Using that $ c(y\ot y^2) = (y^2-x^2)\ot y - \dfrac{1}{2}(1+\xi)(xy + yx)\ot x $, we have that
\begin{align*}
\com(y^3) &= y^3\ot 1 -\dfrac{1}{2}(1+\xi) xy\ot x - y^2\ot y -\dfrac{1}{2}(1-\xi) x\ot xy + y\ot y^2 + 1\ot y^3,
\end{align*}
and consequently $ \com(y^4) = y^4\ot 1 + 1\ot y^4 $,
because $ c(y\ot y^3) = (yx^2 +xyx - y^3 +x^2y)\ot y + \dfrac{1}{2}(1+\xi)(yxy -x^3 + xy^2 + y^2x)\ot x $ 
and $ c(y\ot xy) = \xi xy\ot y + \dfrac{1}{2}(1-\xi)x^2\ot x $. 
Hence, the relation $ y^4 = 0 $ also holds in $\toba(V_{3,1}) $. 
\epf

\begin{prop}\label{prop:BV33}
$\toba(V_{3, 3})$  is the algebra generated by the elements $x, y$ 
satisfying the following relations: 
$$x^2 -2y^2 = 0,\qquad xy + yx = 0,\qquad y^4 = 0.$$
In particular,
$\dim \toba(V_{3, 3}) = 8 $.
\end{prop}
\pf 
Using the braiding given in Proposition \ref{prop:braidedspacesdim2}, we have that
\begin{align*}
\com(x^2) &= x^2\ot 1 +(1-\xi) x\ot x + 1\ot x^2,&
\com(xy) &= xy\ot 1 -\xi x\ot y - y\ot x + 1\ot xy, \\
\com(yx) &= yx\ot 1 + y\ot x +\xi x\ot y + 1\ot yx, &
\com(y^2) &= y^2\ot 1 +\dfrac{1}{2}(1 - \xi) x\ot x + 1\ot y^2.
\end{align*}
This gives us that the relations $ x^2 -2y^2 = 0 $ and $ xy + yx = 0 $
must hold in $\toba(V_{3,3}) $. 
Since $ c(y\ot y^2) = (y^2-x^2)\ot y + \dfrac{1}{2}(\xi -1)(xy + yx)\ot x $, it follows that
\begin{align*}
\com(y^3) &= y^3\ot 1 +\dfrac{1}{2}(\xi -1) xy\ot x - y^2\ot y -
\dfrac{1}{2}(1+\xi) x\ot xy + y\ot y^2 + 1\ot y^3,
\end{align*}
and consequently $ \com(y^4) = y^4\ot 1 + 1\ot y^4 $,
because $ c(y\ot y^3) = (yx^2 +xyx - y^3 +x^2y)\ot y + \dfrac{1}{2}(1-\xi)(yxy -x^3 + xy^2 + y^2x)\ot x $ 
and $ c(y\ot xy) = -\xi xy\ot y + \dfrac{1}{2}(1+\xi)x^2\ot x $. 
\epf

We end this section with the characterization of the finite-dimensional 
Nichols algebras over indecomposable objects 
in $\ydk$.

\subsection*{Proof of Theorem \ref{thm:Nicholskdm}}
Let $V$ be an indecomposable module such that $\toba(V)$ is finite-dimensional. Then
by Theorem \ref{thm:indecnicholsinfdim}, $V$ is necessarily simple. The claim then 
follows by Lemmata 
\ref{lema:nichols-1-dim} and 
\ref{lema:2dim-infinitenichols}, and Propositions \ref{prop:BV21}, \ref{prop:BV23},
\ref{prop:BV31} and \ref{prop:BV33}. Clearly, Nichols algebras over distinct families are pairwise 
non-isomorphic, since they are generated by the set of primitive elements which are
non-isomorphic as Yetter–Drinfeld modules.
\qed

\section{Hopf algebras over $\Kk$}\label{sec:HaK}

In this last section we determine all finite-dimensional Hopf algebras $H$ such that 
$H_{[0]}=\Kk$ and the corresponding infinitesimal braiding is a simple object in $\ydk$ under the assumption that the diagram is a Nichols algebra. 
That is, the graded algebra with respect to the standard filtration is 
$\gr H = \bigoplus_{i\geq 0} H_{[i]}/H_{[i-1]}\simeq \toba (R(1)) \#\Kk$ 
with $R(1)$ isomorphic to a simple object in 
$\ydk$.

Next,
we show that the bosonizations of the Nichols algebras associated with the simple modules 
$\Bbbk_{\chi^{\ell}}$ with $\ell=1,3$ and $V_{2,1}$, $V_{2,3}$ do not admit deformations.

Recall that, for $v \in V= R(1)$, the formula given by the bosonization yields 
\begin{align*}
\com(v \# 1) &= v^{(1)} \# (v^{(2)})_{(-1)} \ot (v^{(2)})_{ (0)} \# 1 =  v \# 1 \ot 1 \# 1 + 1 \# v_{(-1)} \ot v_{ (0)} \# 1.
\end{align*}
We also write $v=v\#1$ for all $v \in V$, and $k=1\#k$ for all $k\in\Kk$.

\begin{prop}\label{prop:nodefdim1}
Let $H$ be a finite-dimensional 
Hopf algebra over $\Kk$ such that its infinitesimal braiding $V$ is isomorphic to 
$\Bbbk_{\chi^{\ell}}$ with $\ell =1$ or $3$. 
Assume that the diagram is a Nichols algebra.
Then $H\simeq (\bigwedge \Bbbk_{\chi^{\ell}})\# \Kk $.
\end{prop}

\pf Write $\bigwedge \Bbbk_{\chi^{\ell}} = \Bbbk[x]/(x^{2})$.
As $\gr H \simeq (\bigwedge \Bbbk_{\chi^{\ell}})\#\Kk$ 
with $\ell = 1$ or $3$, 
we have to prove that the defining relation $ x^2 = 0 $ of $\bigwedge \Bbbk_{\chi^{\ell}}$ remains homogeneous in $H$. 
Since $\delta(x) = a^{2}\ot x$, we obtain that
\begin{align*}
	\com(x^{2}) = (x\ot 1 + a^2\ot x)^2 = x^{2}\ot 1 + 1\ot x^{2} + (a^{2}\cdot x + x) a^2\ot x = 
	x^{2}\ot 1 + 1\ot x^{2},
\end{align*}
which implies that $x^{2}=0$ in $ H $, because $\mathcal{P}(\Kk) = \{0\}$.
\epf

\begin{prop}\label{prop:nodef21-23}
Let $H$ be a finite-dimensional 
Hopf algebra over $\Kk$ such that its infinitesimal braiding $V$ is isomorphic either to 
$V_{2,1}$ or $V_{2,3}$. 
Assume that the diagram is a Nichols algebra. 
Then $H\simeq \mathfrak{B}(V)\# \Kk $.
\end{prop}

\pf We know that  
$\gr H \simeq \toba(V) \#\Kk$ 
with $V$ isomorphic either to 
$V_{2,1}$ or $V_{2,3}$. As in the proof of Proposition \ref{prop:nodefdim1},
we prove that the homogeneous relations also hold in $H$.  

Assume first that $V\simeq V_{2,1}$.
Then, $ \mathfrak{B}(V_{2, 1})\# \Kk $ is the algebra generated by the elements $x, y, a, b, c, d$ with 
$x,y$ satisfying
the relations of $ \mathfrak{B}(V_{2, 1}) $ \eqref{eq:relBV21}, 
 $a, b, c, d$ satisfying the relations of $\Kk$ \eqref{eq:defining_relations_K}, and  all together satisfying the following relations:
\begin{equation}\label{eq:relV21}
\begin{tabular}{cccc}
$ax = -xa,$ &$ay = \xi ya + xc,$ &$bx = -xb,$ &$by = \xi yb + xd,$ \\
$cx = -xc,$ &$cy = -\xi yc + xa,$ &$dx = -xd,$ &$dy = -\xi yd + xb$.
\end{tabular}
\end{equation}
As $ \Delta(x) = x\ot 1+ d\ot x + (\xi -1)c\ot y $ and $
\Delta(y) = y\ot 1+ a\ot y - \dfrac{\xi +1}{2}b\ot x, $
we have that $ \com(xy +\xi yx)  = (xy +\xi yx)\ot 1+ 1\ot (xy +\xi yx) $ and $
\Delta(y^4)  = y^4\ot 1+ 1\ot y^4. $
Since $\mathcal{P}(H) = \{0\}$, it follows that the relations 
$xy +\xi yx = 0$ and $y^{4}=0$ hold in $H$.

On the other hand, $ \com(x^2) = x^2\ot 1+ a^2\ot x^2 +(\xi -1)ab\ot (xy +\xi yx)= x^2\ot 1+ a^2\ot x^2 $, that is, $x^{2}$ is a $(1,a^{2})$-primitive element in $H_{[1]}$.
As $P_{1, a^2}(H_{[1]}) = P_{1, a^2}(\Kk) =\Bbbk\{1-a^2, ab\} $, we must have that 
\[ x^2 = \mu_1 (1-a^2)+ \mu_2 ab\qquad\text{ for some }\mu_{1}, \mu_{2}\in \Bbbk.
\]
But, by \eqref{eq:relV21},
\begin{align*}
0 &= ax^2- x^2a = a(\mu_1 (1-a^2)+ \mu_2 ab) - (\mu_1 (1-a^2)+ \mu_2 ab)a = \mu_2(1+\xi)c, \\
0 &= bx^2- x^2b= b(\mu_1 (1-a^2)+ \mu_2 ab) - (\mu_1 (1-a^2)+ \mu_2 ab)b = 2\mu_1 c,
\end{align*}
which implies that $ \mu_1 = \mu_2 = 0 $. Therefore, the relation $ x^2 =0 $ also holds in $H$ and 
consequently, $H\simeq \gr H$.

For $V\simeq V_{2,3}$, the proof follows the same lines as for $V\simeq V_{2,1}$.
\epf

Hereafter, we define two Hopf
algebras $\A_{3,1}(\mu)$ and $\A_{3,3}(\mu)$, which are constructed by deforming 
the relations on the Nichols algebras $\toba(V_{3,1})$ and $\toba(V_{3,3})$ over $\Kk$,
respectively, and show
that they are liftings of the corresponding bosonizations.

\begin{defi}\label{def:A31} For $\mu \in \Bbbk$, let $\A_{3,1}(\mu)$ be the 
algebra generated by the elements $x, y, a, b, c$, $d$ satisfying the relations \eqref{eq:defining_relations_K} and the following ones:
\begin{align}
\notag ax & = -\xi xa, & ay & = - ya - xc, & bx &= -\xi xb, & by &= - yb - xd, && \\ \label{eqn:relations:A31}
cx & = \xi xc, & cy &= - yc + xa, & dx &= \xi xd, & dy &= - yd + xb,&&
\end{align}
\vspace{-0.7cm}
\begin{align*}
&& x^2 -2 y^2 & = \mu (1-a^2),  &xy+yx &= \xi\mu ac, &
y^4 &= -\mu y^{2}(1-a^2)  - \dfrac{\mu^{2}}{2}(1-a^2).
 \end{align*}
\end{defi}

$\A_{3,1}(\mu)$ is a Hopf algebra with coalgebra structure and antipode  determined by \eqref{eq:defining_coproduct_antipode_K} and:
\begin{align*}
\com(x) & = x\ot 1 + d\otimes x +(\xi-1) c\otimes y, &
\com(y) & = y \ot 1 + a\otimes y +\dfrac{1}{2}( 
-\xi-1) b\otimes x,
\end{align*}
\vspace{-0.7cm}
\begin{align*}
\Ss(x) &= -ax - (1+\xi)cy, &  \Ss(y) &= -dy + \frac{1}{2}(\xi-1)bx, &  \varepsilon(x) &= \varepsilon(y) = 0.
\end{align*}

\begin{obs}\label{rmk:A31} Clearly, $\A_{3,1}(0)\simeq \toba(V_{3,1})\#\Kk$.  
Also note that $\A_{3,1}(\mu)$ is the quotient of the algebra 
$ T(V_{3,1}) \ot \Kk$ by the two-sided ideal generated by the relations \eqref{eqn:relations:A31}; denote this ideal by $ J_{3, 1} $. Furthermore,  
formulae \eqref{eq:rest_of_coproduct_and_antipode_in_K} hold.
\end{obs}

\begin{defi}\label{def:A33}
For $\mu \in \Bbbk$, let $\A_{3,3}(\mu)$ be the 
Hopf algebra defined by $\A_{3,1}(\mu) = \A_{3,3}(\mu)$ 
as algebra but with its coalgebra structure determined by 
the same counit and comultiplication 
for the generators $a,b,c,d$, but 
\begin{align*}
\com(x) & = x\ot 1+ a\ot x + (\xi +1)b\ot y, & \varepsilon(x) & =0,\\
\com(y) & = y\ot 1+ d\ot y + \dfrac{1}{2}(1-\xi)c\ot x,
&\varepsilon(y) &= 0.
\end{align*}
In particular, we have that
\begin{align*}
\Ss(x) &= -dx - (\xi-1)by, &  \Ss(y) &= -ay + \frac{1}{2}(1+\xi)cx.
\end{align*}
\end{defi}

\begin{obs} As before, $\A_{3,3}(0)\simeq \toba(V_{3,3})\#\Kk$ and   $\A_{3,3}(\mu)$ 
is the quotient of the algebra 
$ T(V_{3,3}) \ot \Kk$ by the two-sided ideal generated by the relations \eqref{eqn:relations:A31}.  
\end{obs}

In the next lemma we show that the algebras $\A_{3,1}(\mu)$ and $\A_{3,3}(\mu)$
are finite-dimensional Hopf algebras over $\Kk$.

\begin{lema}\label{lema:A31-A33overK}  
Let $\mu \in \Bbbk$. The Hopf algebras $\A_{3,1}(\mu)$ and $\A_{3,3}(\mu)$ are 
finite-dimensional and $(\A_{3,1}(\mu))_{[0]} \simeq \Kk\simeq (\A_{3,3}(\mu))_{[0]} $.
\end{lema}

\pf We prove the assertion for $\A_{3,1}(\mu) $, being the proof for $\A_{3,3}(\mu)$ 
completely analogous. First note that, by Remark \ref{rmk:A31}, $\A_{3,1}(\mu) = T(V_{3,1})\ot \Kk/J_{3,1}$. 
Also note that  $T(V_{3,1})\ot \Kk$ is a graded algebra with the gradation defined as usual on $T(V_{3,1})$
and all the elements in $\Kk$ being of degree $0$.

Let $\A_{0}$ be the subalgebra generated by the coalgebra $C$ linearly spanned by $a,b,c,d$. Then
$\A_{0}$ is a Hopf subalgebra of $\A_{3,1}(\mu)$. We claim that $\A_{0}$ is isomorphic to $\Kk$. Indeed,
consider the Hopf algebra map $\varphi: \Kk \to \A_{3,1}(\mu)$ given by the composition 
$\Kk \hookrightarrow T(V_{3,1})\ot \Kk \twoheadrightarrow T(V_{3,1})\ot \Kk/J_{3,1} $. Clearly, 
$\A_{0} \simeq \varphi (\Kk)$. Since 
$\dim \Kk = 8$, to prove that $\varphi(\Kk) \simeq \Kk$ it is enough to show that 
$\dim \varphi(\Kk) > 4$. 
As the relations in $J_{3,1}$ do not involve relations in degree $0$, it follows that $\varphi(C) \simeq C$. 
Moreover, since the elements $a,b,c,d$ are linearly independent in $\Kk$, they are also so in 
$\A_{3,1}(\mu)$. Hence, $\varphi(\Kk) \simeq  \Kk$ and the claim is proved.

Set
$\A_{1} = \A_{0} +\Kk \{x,y\}$, 
$\A_{2} = \A_{1} + \Kk\{xy, y^{2}\}$, $\A_{3}=\A_{2}+\Kk\{xy^{2}, y^{3}\}$ and 
$\A_{4}=\A_{3} + \Kk\{xy^{3}\}$. Then 
$\{\A_{n}\}_{0\leq n\leq 3}$ is a coalgebra filtration of $\A_{3,1}(\mu)$. In particular, 
$(\A_{3,1}(\mu))_{0} \subseteq \Kk$ and consequently $(\A_{3,3}(\mu))_{[0]} = \Kk$;
that is, $\A_{3,1}(\mu)$ is a Hopf algebra over $\Kk$. Hence, $\A_{3,1}(\mu)$ is a finite-dimensional
Hopf algebra which is free over $\Kk$. So, $8$ divides $\dim \A_{3,1}(\mu)$.
Besides, $\A_{3,1}$ is a $\Kk$-module with 
a set of generators $\{1,x,y, xy,y^{2}, xy^{2}, y^{3}, xy^{3}\}$. Thus, 
$\dim \A_{3,1}(\mu) \leq 8\dim \Kk = 64$ and the lemma is proved.
\epf

Now, we show that the algebras $\A_{3,1}(\mu)$ and $\A_{3,3}(\mu)$
are liftings of 
$\toba(V_{3,1})\# \Kk$ and 
$\toba(V_{3,3})\# \Kk$ for all $\mu \in \Bbbk$, respectively. 

\begin{lema} \label{lema:A31-A33lioftingsoverK}
$\gr \A_{3,1}(\mu) \simeq \toba(V_{3,1})\# \Kk$ and $\gr \A_{3,3}(\mu) \simeq \toba(V_{3,3})\# \Kk$.
\end{lema}

\pf It is enough to show that $\dim \A_{3,1}(\mu), \dim \A_{3,3}(\mu) \geq 64$,
since by the proof of Lemma \ref{lema:A31-A33overK} we have that 
$\gr \A_{3,1}(\mu) \simeq R_{3,1}\# \Kk $ and $\gr \A_{3,3}(\mu) \simeq R_{3,3}\# \Kk $,
where $R_{3,1}$ and $R_{3,3}$ are $K$-modules linearly spanned by the set 
$\{1,x,y, xy,y^{2}, xy^{2}, y^{3}, xy^{3}\}$. 
We show that the set $B= \{x^iy^ja^kb^l: 0\leq i, l\leq 1,0\leq j, k\leq 3\}$ 
is linearly 
independent by using appropriate representations.
As $ \A_{3,1}(\mu) $ and $ \A_{3,3}(\mu) $ have the same algebra structure, we prove it
only for $ \A_{3,1}(\mu) $. 

For $ \lambda $ a $4$-th root of unity, consider the 
$8$-dimensional representation $W_{\lambda }$ of $ \A_{3,1}(\mu) $ determined by the following matrices on a fixed basis of 
$W_{\lambda }$:
\begin{align*}
\rho_1(a) & = \left(
\begin{smallmatrix}
\textbf{a} & \textbf{0} \\
\textbf{0} & -\xi \textbf{a} \\
\end{smallmatrix}
\right), & 
\rho_1(b) & = \left(
\begin{smallmatrix}
\textbf{0} & \textbf{b}_{12} \\
\textbf{b}_{21} & \textbf{0} \\
\end{smallmatrix}
\right), &
\rho_1(x) & = \left(
\begin{smallmatrix}
\textbf{0} & \textbf{x} \\ 
\id_4 & \textbf{0} \\ 
\end{smallmatrix}
\right), &
\rho_1(y) & = \left(
\begin{smallmatrix}
\textbf{y}_{11} & \textbf{0} \\ 
\textbf{0} & \textbf{y}_{22} \\ 
\end{smallmatrix}
\right),
\end{align*}
where 
\begin{align*}
\textbf{a} &= \left(
\begin{smallmatrix}
\lambda & 0 & \mu(\lambda^3-\lambda) & 0 \\
0 & -\lambda & 0 & \mu(\lambda-\lambda^3) \\
0 & 0 & -\lambda & 0  \\
0 & 0 & 0 & \lambda  \\
 \end{smallmatrix}
\right),&
\textbf{b}_{12}&= \left(
\begin{smallmatrix}
0 & \xi\mu(\lambda^3-\lambda) & 0 & \xi\mu^2(\lambda-\lambda^3) \\
0 & 0 & 0 & 0 \\
0 & 2\xi \lambda^3 & 0 & \xi\mu(\lambda-\lambda^3) \\
0 & 0 & 0 & 0 \\
 \end{smallmatrix}
\right),\\
\textbf{b}_{21}&= \left(
\begin{smallmatrix}
0 & -\lambda^3 & 0 & \mu(\lambda^3-\lambda) \\
0 & 0 & 0 & 0  \\
0 & 0 & 0 & \lambda^3 \\
0 & 0 & 0 & 0  \\
 \end{smallmatrix}
\right),&
\textbf{x}&= \left(
\begin{smallmatrix}
\mu(1-\lambda^2) & 0 & \mu^2(\lambda^2-1) & 0 \\ 
0 & \mu(1-\lambda^2) & 0 & \mu^2(\lambda^2-1) \\ 
2 & 0 & \mu(\lambda^2-1) & 0 \\ 
0 & 2 & 0 & \mu(\lambda^2-1) \\ 
 \end{smallmatrix}
\right),\\
\textbf{y}_{11}&= \left(
\begin{smallmatrix}
0 & 0 & 0 & \frac{1}{2}\mu^2(\lambda^2-1) \\ 
1 & 0 & 0 & 0 \\ 
0 & 1 & 0 & \mu(\lambda^2-1) \\ 
0 & 0 & 1 & 0 \\
 \end{smallmatrix}
\right),&
\textbf{y}_{22}&= \left(
\begin{smallmatrix}
0 & \mu\lambda^2 & 0 & \frac{1}{2}\mu^2(1-\lambda^2) \\ 
-1 & 0 & 0 & 0 \\ 
0 & -1 & 0 & \mu \\ 
0 & 0 & -1 & 0 \\
 \end{smallmatrix}
\right).
\end{align*}

Assume that
\begin{equation}\label{linear combination}
\sum\limits_{0\leq i, l\leq 1,0\leq j, k\leq 3}f_{i,j,k,l} x^iy^ja^kb^l = 0,
\end{equation} 
for some $f_{i,j,k,l} \in \Bbbk$. Applying \eqref{linear combination} to the first vector of the fixed basis, we get that:
\[ f_{i,j,0,0} +\lambda f_{i,j,1,0}+\lambda^2 f_{i,j,2,0}+\lambda^3 f_{i,j,3,0} = 0, 
\ \text{ for all }\ 0\leq i\leq 1,0\leq j\leq 3.
\]
Since this equation must hold for any $4$-th root of unity $ \lambda $, it follows that $ f_{i,j,k,0} = 0 $ 
for all $ 0\leq i\leq 1$ and $0\leq j, k\leq 3 $. 
To prove that the remaining coefficients are zero, 
we need another representation. 

For $ \lambda $ a $4$-th root of unity, consider
now the $16$-dimensional representation 
$U_{\lambda}$ given by the following matrices on a fixed basis of $U_{\lambda}$:
\begin{align*}
\rho_2(a) &= \left(
\begin{smallmatrix}
\textbf{a} & \textbf{0} & \mu(\textbf{d}-\textbf{a}) & \textbf{0} & \textbf{0} & \xi\mu(\textbf{c}-\textbf{b}) & \textbf{0} & \textbf{0} \\
\textbf{0} & -\textbf{a} & \textbf{0} & \mu(\textbf{a}-\textbf{d}) & \textbf{0} & \textbf{0} & \textbf{0} & \textbf{0} \\
\textbf{0} & \textbf{0} & -\textbf{a} & \textbf{0} & \textbf{0} & 2\xi \textbf{c} & \textbf{0} & -\xi\mu (\textbf{b}+\textbf{c}) \\
\textbf{0} & \textbf{0} & \textbf{0} & \textbf{a} & \textbf{0} & \textbf{0} & \textbf{0} & \textbf{0} \\
\textbf{0} & -\textbf{c} & \textbf{0} & \mu \textbf{c} & -\xi \textbf{a} & \textbf{0} & \xi\mu(\textbf{a}-\textbf{d}) & \textbf{0} \\
\textbf{0} & \textbf{0} & \textbf{0} & \textbf{0} & \textbf{0} & \xi \textbf{a} & \textbf{0} & \xi\mu(\textbf{d}-\textbf{a})  \\
\textbf{0} & \textbf{0} & \textbf{0} & \textbf{c} & \textbf{0} & \textbf{0} & \xi \textbf{a} & \textbf{0}  \\
\textbf{0} & \textbf{0} & \textbf{0} & \textbf{0} & \textbf{0} & \textbf{0} & \textbf{0} & -\xi \textbf{a}  \\
\end{smallmatrix}
\right),\\
\rho_2(b) &= \left(
\begin{smallmatrix}
\textbf{b} & \textbf{0} & -\mu \textbf{b} & \textbf{0} & \textbf{0} & \xi\mu(\textbf{d}-\textbf{a}) & \textbf{0} & \xi\mu^2(\textbf{a}-\textbf{d}) \\
\textbf{0} & -\textbf{b} & \textbf{0} & \mu \textbf{b} & \textbf{0} & \textbf{0} & \textbf{0} & \textbf{0} \\
\textbf{0} & \textbf{0} & -\textbf{b} & \textbf{0} & \textbf{0} & 2\xi \textbf{d} & \textbf{0} & \xi\mu (\textbf{a}-\textbf{d}) \\
\textbf{0} & \textbf{0} & \textbf{0} & \textbf{b} & \textbf{0} & \textbf{0} & \textbf{0} & \textbf{0} \\
\textbf{0} & -\textbf{d} & \textbf{0} & \mu (\textbf{d}-\textbf{a}) & -\xi \textbf{b} & \textbf{0} & \xi\mu \textbf{b} & \textbf{0} \\
\textbf{0} & \textbf{0} & \textbf{0} & \textbf{0} & \textbf{0} & \xi \textbf{b} & \textbf{0} & -\xi\mu \textbf{b}  \\
\textbf{0} & \textbf{0} & \textbf{0} & \textbf{d} & \textbf{0} & \textbf{0} & \xi \textbf{b} & \textbf{0}  \\
\textbf{0} & \textbf{0} & \textbf{0} & \textbf{0} & \textbf{0} & \textbf{0} & \textbf{0} & -\xi \textbf{b}  \\
\end{smallmatrix}
\right), &
\rho_2(x)&= \left(
\begin{smallmatrix}
\textbf{0}_8 & \textbf{x} \\
\id_8 & \textbf{0}_8  \\
\end{smallmatrix}
\right),\\
\rho_2(y)&= \left(
\begin{smallmatrix}
\textbf{0} & \textbf{0} & \textbf{0} & \frac{\mu^2}{2}(\textbf{a}^2-\id_2) & \xi\mu \textbf{ac} & \textbf{0} & -\xi\mu^2 \textbf{ab} & \textbf{0} \\
\id_2 & \textbf{0} & \textbf{0} & \textbf{0} & \textbf{0} & \xi\mu \textbf{ac} & \textbf{0} & -\xi\mu^2 \textbf{ab} \\
\textbf{0} & \id_2 & \textbf{0} & \mu(\textbf{a}^2-\id_2) & \textbf{0} & \textbf{0} & \xi\mu \textbf{ac} & \textbf{0} \\
\textbf{0} & \textbf{0} & \id_2 & \textbf{0} & \textbf{0} & \textbf{0} & \textbf{0} & \xi\mu \textbf{ac} \\
\textbf{0} & \textbf{0} & \textbf{0} & \textbf{0} & \textbf{0} & \mu \textbf{a}^2 & \textbf{0} & \frac{\mu^2}{2}(\id_2-\textbf{a}^2) \\
\textbf{0} & \textbf{0} & \textbf{0} & \textbf{0} & -\id_2 & \textbf{0} & \textbf{0} & \textbf{0}  \\
\textbf{0} & \textbf{0} & \textbf{0} & \textbf{0} & \textbf{0} & -\id_2 & \textbf{0} & \mu \id_2  \\
\textbf{0} & \textbf{0} & \textbf{0} & \textbf{0} & \textbf{0} & \textbf{0} & -\id_2 & \textbf{0}  \\
\end{smallmatrix}
\right),
\end{align*}
where 
$
\textbf{a} = \left(
\begin{smallmatrix}
\lambda & 0 \\
0 & -\xi\lambda  \\
\end{smallmatrix}
\right), \ 
\textbf{b}= \left(
\begin{smallmatrix}
0 & \lambda^2 \\
0 & 0 \\
\end{smallmatrix}
\right), \ 
\textbf{c}= \left(
\begin{smallmatrix}
0 & 1 \\
0 & 0 \\
\end{smallmatrix}
\right), \
\textbf{d} = \left(
\begin{smallmatrix}
\lambda^3 & 0 \\
0 & \xi\lambda^3  \\
\end{smallmatrix}
\right)$ and
\[
\textbf{x}= \left(
\begin{smallmatrix}
\mu(\id_2-\textbf{a}^2) & \textbf{0} & \mu^2(\textbf{a}^2-\id_2) & \textbf{0} \\
\textbf{0} & \mu(\id_2-\textbf{a}^2) & \textbf{0} & \mu^2(\textbf{a}^2-\id_2) \\
2\id_2 & \textbf{0} & \mu(\textbf{a}^2-\id_2) & \textbf{0} \\
\textbf{0} & 2\id_2 & \textbf{0} & \mu(\textbf{a}^2-\id_2) \\
\end{smallmatrix}
\right).\]
Applying the residual equation of (\ref{linear combination}) to the second vector of the 
fixed basis, we get that
\[ \lambda^2(f_{i,j,0,1} +\lambda f_{i,j,1,1}+\lambda^2 f_{i,j,2,1}+\lambda^3 f_{i,j,3,1}) = 0, 
\ \text{ for all }\  0\leq i\leq 1,0\leq j\leq 3,
\]
implying that $ f_{i,j,k,1} = 0 $ for all $ 0\leq i\leq 1,0\leq j, k\leq 3 $. 
Therefore, $B$ is a linearly independent set and $\A_{3,1}(\mu)$ is a lifting of $ \toba(V_{3,1})\# \Kk $,
for all $\mu \in \Bbbk$.
\epf

\begin{prop}\label{prop:def31-33}
Let $H$ be a finite-dimensional Hopf algebra over $\Kk$ such that its infinitesimal braiding is
isomorphic to $V_{3,1}$ or $V_{3,3}$. 
Assume that the diagram is a Nichols algebra. 
Then $H \simeq \A_{3,1}(\mu)$ or $H \simeq \A_{3,3}(\mu)$
for some $\mu \in \Bbbk$, respectively.
\end{prop}

\pf We have that $\gr H \simeq \toba(V)\# \Kk$ with 
$V\simeq V_{3,1}$ or $V\simeq V_{3,3}$.
Recall that $ \mathfrak{B}(V)\# \Kk $ is the algebra generated by $x, y, a, b, c, d$, where 
$x,y$ are the 
generators of $ \mathfrak{B}(V) $, 
$a, b, c, d$ are the generators of 
$\Kk$, and they all together satisfy the first two rows of relations \eqref{eqn:relations:A31}.

We prove the claim for $V\simeq V_{3,1}$. The proof for $V\simeq V_{3,3}$ follows the same lines. 
As $\Delta(x) = x\ot 1+ d\ot x + (\xi -1)c\ot y $ and $\Delta(y) = y\ot 1+ a\ot y - \dfrac{\xi +1}{2}b\ot x$,
we obtain that
\begin{align}\label{eq:coproducts31_1}
\Delta(x^2-2y^2) &= (x^2-2y^2)\ot 1+ a^2\ot (x^2-2y^2) \qquad \text{ and }\\
\label{eq:coproducts31_2}\Delta(xy+yx) &= (xy+yx)\ot 1+ 1\ot (xy+yx) -\xi ac\ot (x^2-2y^2).
\end{align}
By \eqref{eq:coproducts31_1}, we get that 
$ x^2 -2 y^2\in P_{1, a^2}(H_{[1]}) 
= P_{1, a^2}(\Kk) =\Bbbk\{1-a^2, ab\} $. Then, in $ H $, 
\[ x^2 -2 y^2 = \mu_1 (1-a^2)+ \mu_2 ab \qquad\text{ for some }\mu_{1}, \mu_{2}\in \Bbbk.
\]
Thus, \eqref{eq:coproducts31_2} can be rewritten as
\begin{align*}
\com(xy+yx-\xi \mu_{1} ac) &= (xy+yx - \mu_{1}\xi ac)\ot 1+ 1\ot (xy+yx- \mu_{1}\xi ac) -\xi\mu_2 ac\ot ab.
\end{align*}
However, a tedius calculation on $ H_{[1]} $ shows that $ \mu_2 $ must be $ 0 $ in the last equation. Hence,
\[ xy+yx = \xi\mu_1 ac \qquad \text{ and } \qquad x^2 -2 y^2 = \mu_1 (1-a^2) \qquad\text{ for some }\mu_{1}\in \Bbbk.
\]
Finally, observe that the element $ R := y^4 +\mu_1 y^2(1-a^2) + \dfrac{1}{2}\mu_1^2(1-a^2)$ satisfies $ \Delta(R) = R\ot 1 + 1\ot R $, which implies that $ R= 0 $ in $ H $.

Since the defining relations of $\A_{3,1}(\mu_{1})$
hold in $H$,
there exists a surjective
Hopf algebra map from $H$ to $\A_{3,1}(\mu_{1})$. As both algebras have the same dimension, 
the proposition follows.
\epf

We end the paper with the proof of our second main theorem.

\subsection*{Proof of Theorem \ref{mthm:liftingsoverK}} 
By Theorem \ref{thm:Nicholskdm}, we have that $\gr H \simeq \toba(V)\# \Kk$ with $V$ isomorphic to $\Bbbk_{\chi}$, $\Bbbk_{\chi^{3}}$, $V_{2,1}$, $V_{2,3}$, $V_{3,1}$ or $V_{3,3}$.

If $V \simeq \Bbbk_{\chi}$, $\Bbbk_{\chi^{3}}$, $V_{2,1}$ or $V_{2,3}$, then 
$H\simeq \toba(V)\# \Kk$ by Propositions \ref{prop:nodefdim1} and \ref{prop:nodef21-23}.
If $V\simeq V_{3,1}$ or $V\simeq V_{3,3}$, then by Proposition \ref{prop:def31-33} it follows
that $H\simeq \A_{3,1}(\mu)$ or 
$H\simeq \A_{3,3}(\mu)$ for some $\mu \in \Bbbk$, respectively.

Conversely, it is clear that the algebras listed in items $(i)$, $(ii)$ and $(iii)$ are liftings
of Hopf algebras over $\Kk$. The Hopf algebras $\A_{3,1}(\mu)$ and 
$\A_{3,3}(\mu)$ are also liftings by
Lemma \ref{lema:A31-A33lioftingsoverK}.

Finally, two algebras from different families are not isomorphic as Hopf
algebras since their infinitesimal braidings are not isomorphic as Yetter–Drinfeld
modules.
\qed

\section*{Acknowledgements}
The authors thank
L. Vendramin for providing the computation with GAP of 
the Nichols algebras over $2$-dimensional modules, and 
A. Garc\'ia Iglesias, N. Andruskiewitsch and I. Angiono for 
fruitful discussions on the results of the paper. The authors also wish to thank
the referee for his/her valuable comments that helped to improve considerably the paper.
Part of the proofs of Remark \ref{rmk:soctopone-dim} and Lemma \ref{lem:3dimindec} are due to him/her.

\end{document}